\documentclass[a4paper,11pt]{article}
\usepackage{geometry}
\usepackage{mathtext}
\usepackage[latin2]{inputenc}
\usepackage[english]{babel}
\usepackage{bbm}
\usepackage{dsfont}
\usepackage{tikz-cd}
\usepackage{amsmath}
\usepackage{enumitem}
\usepackage{amsfonts}
\usepackage{amssymb}
\usepackage{mathrsfs}
\usepackage{amsthm}
\usepackage{euscript}
\usepackage{authblk}
\usepackage{stmaryrd}
\usepackage{graphicx}
\usepackage{scalerel}
\usepackage{wasysym}
\usepackage{enumitem}

\usepackage[format=plain,
            labelfont={bf,it},
            textfont=it]{caption}

\usepackage [english]{babel}
\usepackage [autostyle, english = american]{csquotes}
\MakeOuterQuote{"}

\DeclareMathOperator{\supp}{supp}

\DeclareMathOperator{\RR}{\mathbb{R}}
\DeclareMathOperator{\NN}{\mathbb{N}}

\DeclareMathOperator{\Pp}{\mathbb{P}}
\DeclareMathOperator{\EE}{\mathbbm{E}}

\DeclareMathOperator{\chs}{\mathds{1}}
\DeclareMathOperator{\AAA}{\mathscr{A}}
\DeclareMathOperator{\FF}{\mathcal{F}}

\DeclareMathOperator{\cgc}{\mathcal{G}}
\DeclareMathOperator{\Ec}{\mathcal{E}}
\DeclareMathOperator{\Tc}{\mathcal{T}}
\DeclareMathOperator{\Ucc}{\mathcal{U}}

\DeclareMathOperator{\Rc}{\mathcal{R}}

\makeatletter
\newcommand{\precline}{\mathrel{\mathpalette\pr@ceqd@t\relax}}
\newcommand{\pr@ceqd@t}[2]{%
  \begingroup
  \sbox\z@{$#1\prec$}\sbox\tw@{$#1\prec$}%
  \dimen@=\dimexpr\ht\tw@-\ht\z@\relax
  {\prec}%
  \mkern-5mu
  \raisebox{\dimen@}{$\m@th#1\; |$}%
  \endgroup
}
\makeatother
\newcommand{\khan}[1]{H_{#1}}

\newcommand{\DD}[1]{\mathrm{D}_{#1}}

\newcommand{\col}[1]{\mathcal{B}_{#1}}

\newcommand{\eLpx}[1]{\left\|{#1}\right\|_{L^p(\Omega, X)}}
\newcommand{\eLpxb}[1]{\bigg\|{#1}\bigg\|_{L^p(\Omega, X)}}

\newcommand{\eLpxA}[2]{\left\|{#1}\right\|_{L^p(#2, X)}}

\newcommand{\eLp}[1]{\left\|{#1}\right\|_{L^p(\Omega)}}

\newcommand{\eLpp}[3]{\left\|{#1}\right\|_{L^{#3}(#2)}}
\newcommand{\eLppp}[3]{\left\|{#1}\right\|^{#3}_{L^{#3}(#2)}}
\newcommand{\nx}[1]{\left\|{#1}\right\|_{X}}
\newcommand{\nxb}[1]{\bigg\|{#1}\bigg\|_{X}}

\newcommand{\Nay}[1]{N(#1)}
\DeclareMathOperator{\Io}{A_0}
\DeclareMathOperator{\Wb}{\tilde{B}}
\DeclareMathOperator{\Wv}{\tilde{V}}
\DeclareMathOperator{\Wu}{\tilde{U}}

\newcommand{\haxi}[2]{\hai{#1}(#2)}
\newcommand{\hai}[1]{\phi_{#1}}
\newcommand{\ha}[1]{k_{#1}}
\newcommand{\hax}[2]{\ha{#1}(#2)}
\newcommand{\haar}[1]{h_{#1}}

\newcommand{\st}[1]{{#1}^{\tiny{\diamond}}}

\newcommand{\ti}[1]{{#1}^{\tiny{\sun}}}

\newcommand{\crl}[1]{\llbracket{#1}\rrbracket}
\newcommand{\ccc}{\mathscr{C}}
\newcommand{\dyad}[1]{\mathcal{I}^{#1}}

\newcommand{\bbb}{\mathscr{B}}
\newcommand{\GC}[3]{G_{#1}(#2, #3)}
\newcommand{\Est}[1]{{#1}^*}
\newcommand{\dx}{\operatorname{d}\!x}
\newcommand{\dP}{\operatorname{d}\!\Pp}

\newcommand{\stlemcrl}[1]{\frac{(4\crl{#1}+1)^{1-\frac{1}{p}}}{1-2^{-\frac{1}{p}}}}

\newcommand{\stlemcrlp}[1]{\frac{(4\crl{#1}+1)^{p-1}}{(1-2^{-\frac{1}{p}})^{p}}}

\newcommand{\stlemp}{2^{p-1}(1+MT_p(\RR)^{p})}
\newcommand{\stlem}{2^{1-\frac{1}{p}}(1+MT_p(\RR)^{p})^{\frac{1}{p}}}

\newcommand{\dPx}[1]{\operatorname{d}\!\Pp(#1)}
\renewcommand{\leq}{\leqslant}
\renewcommand{\geq}{\geqslant}

\newcommand{\mardifold}[2]{\mathbf{\Delta}_{#2} #1 }

\newcommand{\mardif}[2]{\EE(#1\big|\FF_{#2})-\EE(#1\big|\FF_{#2-1}) }
\newcommand{\mardifx}[3]{\EE(#1\big|\FF_{#2})(#3)-\EE(#1\big|\FF_{#2-1})(#3) }
\newcommand{\mardifdyad}[2]{\EE(#1\big|\dyadsigma{#2})-\EE(#1\big|\dyadsigma{#2-1}) }

\newcommand{\dyadsigma}[1]{\mathcal{D}_{#1}}

 	\definecolor{cornflowerblue}{rgb}{0.39, 0.58, 0.93}
\definecolor{bostonuniversityred}{rgb}{0.8, 0.0, 0.0}
 	\definecolor{skymagenta}{rgb}{0.81, 0.44, 0.69}
   	\definecolor{pumpkin}{rgb}{1.0, 0.46, 0.09}
    \definecolor{deepmagenta}{rgb}{0.8, 0.0, 0.8}

\newcommand{\chng}[1]{\textcolor{black}{#1}}

\theoremstyle{plain}\newtheorem{T1}{Theorem}
\theoremstyle{plain}
\theoremstyle{plain}\newtheorem{L1}[T1]{Lemma}
\theoremstyle{plain}
\theoremstyle{definition}\newtheorem{Def}[T1]{Definition}
\theoremstyle{plain}
\theoremstyle{plain}
\theoremstyle{definition}\newtheorem{Rem}[T1]{Remark}

\title{Martingale Type, the Gamlen-Gaudet Construction and a Greedy Algorithm }
\author[1]{Krystian Kazaniecki and Paul F.X. M\"{u}ller}

\definecolor{grg}{rgb}{0.18,0.49,0.47}

\begin{document}

\maketitle
\begin{abstract}
    In the present paper we identify those filtered probability spaces $(\Omega,\, \FF,\, \left(\FF_n\right),\, \Pp)$ that determine already the martingale type of \chng{a} Banach space $X$. We isolate intrinsic conditions on the filtration $(\FF_n)$ of purely atomic $\sigma$-algebras which determine that the upper $\ell^p$ estimates
    \[
 \eLpp{f}{\Omega,\, X}{p}^p\leq C^p\left( \chng{\eLpp{\EE(f|\FF_0)}{\Omega,\, X}{p}^p}+\sum_{n=1}^{\infty}  \eLpp{\mardif{f}{n}}{\Omega,\, X}{p}^p\right),\qquad f\in L^p(\Omega,X)
\]
imply that the Banach space X is of \chng{martingale} type $p$. Our paper complements \mbox{G. Pisier's} investigation \cite{Pisier1975} and continues the work by S. Geiss and second named author in  \cite{Geiss2008}.
\end{abstract}
\vskip5mm
\let\thefootnote\relax\footnotetext{This research was partially supported by the National Science Centre, Poland, and Austrian Science Foundation FWF joint CEUS programme. National Science Centre project no. 2020/02/Y/ST1/00072 and FWF project no. I5231.

}
\section{Introduction}
Let $(\Omega,\, \FF,\, \Pp)$ be a probability space.
Let $(\FF_n)$ be an increasing sequence of  purely atomic sub-$\sigma$-algebras $\FF$. If $\FF$ is the smallest $\sigma$-algebra containing $\bigcup \FF_n$, then we say that $(\Omega,\, \FF,\, \left(\FF_n\right),\, \Pp)$ forms a filtered probability space. Given a Banach space $X$ we let $L^{p}(\Omega,\,\FF,\,\Pp\,; X)$ denote the Banach space of $p$-integrable Bochner measurable functions.
If the underlying filtration and measure are clear from the context we abbreviate our notation to $L^p(\Omega,X)$ or $L^p(X)$ and $L^p(\Omega)$ in case $X=\RR$. 

By way of introduction we recall Burkholder's classical inequality for scalar valued martingales. If $1<p<\infty$ there exist $c_p>0$, and $C_p<\infty$ such that
\begin{equation}\label{burkholder}
c_p\eLp{f}\leq \eLp{\left(|\EE f|^2+\sum_{n=1}^{\infty}\big|\mardifold{f}{n}\big|^2\right)^{\frac{1}{2}}}\leq C_p \eLp{f},    
\end{equation}
where $\mardifold{f}{n}=\mardif{f}{n}$. The estimate \eqref{burkholder} implies that for $1<p<\infty$ martingale differences converge unconditionally in $L^p(\Omega,\,\FF,\, \Pp)$ i.e.
\begin{equation}\label{uncoditional}
\eLp{\sum_{n=1}^{\infty}\varepsilon_n \mardifold{f}{n}}\leq C_p \eLp{f}.
\end{equation}
for $f\in L^p(\Omega,\,\FF,\, \Pp)$ and $\varepsilon_n\in \{-1,1\}$. Moreover for $1<p\leq 2$ inequality \eqref{burkholder} yields the upper $\ell^p$ estimates in term of martingale differences,
\begin{equation}\label{typescalar}
    \eLp{f}\leq \tilde{C}_p\left(\chng{\eLppp{\EE( f|\FF_0)}{\Omega}{p}}+\sum_{n=1}^{\infty} \eLp{\mardifold{f}{n}}^p\right)^{\frac{1}{p}}.
\end{equation}
Passing from scalar valued martingales to vector valued ones we fix a Banach space $X$. It is well known that neither Burkholder's inequalities nor its consequences \eqref{uncoditional} and \eqref{typescalar} hold true for vector valued martingales in general (see \cite{Pisier2016}). 

The validity of \eqref{uncoditional} respectively \eqref{typescalar} each define severely restricting isomorphic invariants on the underlying Banach space $X$.
 We say that $X$ is a UMD space if there exists $C_p<\infty$ such that for any filtered probability space  $(\Omega,\, \FF,\, \left(\FF_n\right),\, \Pp)$ and any $f\in L^p(X)$
and $\varepsilon_n\in \{-1,1\}$ we have
 \begin{equation}\label{uncoditionalvector}
\eLpx{\sum_{n=1}^{\infty}\varepsilon_n \mardifold{f}{n}}\leq C_p \eLpx{f} 
\end{equation}
where $\mardifold{f}{n}=\mardif{f}{n}$. By $UMD(X,p)$ we denote the smallest constant for which inequality \eqref{uncoditionalvector} is satisfied. As pointed out by Pisier the UMD property of a Banach space $X$ is independent of $1<p<\infty$ (see  \cite{zbMATH03482989, zbMATH03482990}, \cite{Pisier2016}).

A Banach space $X$ satisfies martingale type $p$ if there exists $C_p<\infty$ such that any filtered probability space $(\Omega,\, \FF,\, \left(\FF_n\right),\, \Pp)$ gives rise to the upper $\ell^p$ estimates for $f\in L^p(X)$
\begin{equation}\label{type}
    \eLpx{f}\leq C_p\left( \chng{\eLpp{\EE(f|\FF_0)}{\Omega,\, X}{p}^p}+\sum_{n=1}^{\infty} \eLpx{\mardif{f}{n}}^p\right)^{\frac{1}{p}}.
\end{equation}
 In that case we write $X$ satisfies $MT_p$ and by $MT_p(X)$ we denote smallest constant $C_p$ such that \eqref{type} is satisfied. The inequality \eqref{typescalar} states that $\RR$ is a Banach space of martingale type $p$ for any $p\in (1,2].$ 

 We let $\mathcal{I}_n$ denote the collection of pairwise disjoint dyadic intervals of measure $2^{-n}$ contained in $[0,1)$. Thus $I\subseteq [0,1)$ is a dyadic interval in $\mathcal{I}_n$ if there exists $k\in\{0,1,\ldots, 2^n -1\}$ such that
    \[
    I=\left[\frac{k-1}{2^n},\frac{k}{2^n}\right).
    \]
 We denote by $\mathcal{I}=\bigcup\limits_{n=1}^{\infty}\mathcal{I}_n$. Let $\dyadsigma{n}$ be the $\sigma$-algebra generated by $\mathcal{I}_n$. For a dyadic interval $I\in \mathcal{I}$ we let $I^+\subsetneq I$, respectively $I^{-}$ denote the left, respectively the right, half of $I$. Thus $I^+$, $I^{-}$ are dyadic intervals satisfying $I^{+}\cap I^{-}=\emptyset$ and $I^{+}\cup I^{-}=I$. We define the Haar function $h_I\;:[0,1)\rightarrow \{-1,0,1\}$ by putting
\[
h_I=\chs_{I^+}-\chs_{I^-}.
\]

 We put $\dyadsigma{}=\sigma\left(\bigcup\dyadsigma{n}\right)$ and denote by $([0,1),\,\dyadsigma{},\,(\dyadsigma{n}),\,dt)$ the dyadic filtration of the unit interval. Pisier proved in \cite{Pisier1975} that the dyadic filtration  determines already the martingale type of a Banach space $X$. Indeed in Pisier's article combining Proposition 2.4 b) with Theorem 3.1 b) gives that a Banach space X satisfies martingale type $p$ if and only if there exists $C_p<\infty$ such that
\begin{equation}\label{eq:dyadtype}
 \eLpp{f}{[0,1),\, X}{p}\leq C_p\left(\chng{\eLpp{\EE(f|\dyadsigma{0})}{[0,1),\, X}{p}^p}+\sum_{n=1}^{\infty}  \eLpp{\mardifdyad{f}{n}}{[0,1),\, X}{p}^p\right)^{\frac{1}{p}}
\end{equation}
for $f\in L^p([0,1),X),$ where $L^p([0,1),\,X)=L^p([0,1),\,\dyadsigma{},\,dt,\, X)$. This equivalence is also presented in Theorem 10.22 of Pisier's recent book on martingales in Banach spaces \cite{Pisier2016}. 

We interpret the statement of the above cited Theorem of Pisier in terms of Haar functions. Since $1< p\leq 2$ it is clear that for a given $f\in L^p(X)$ and $n\in\NN$ we have
\[
\mardifdyad{f}{n}=\sum_{I\in \dyad{}_{n-1}} x_I h_I \qquad\mbox{where}\qquad x_I=\int_0^1 f h_I \frac{dx}{|I|},
\]
and consequently 
\[
f=\EE f + \sum\limits_{n=1}^{\infty}\sum_{I\in \dyad{}_{n-1}} x_I h_I
\]
where the series converges in $L^p(X)$. Therefore we may restate \eqref{eq:dyadtype} as
\[
 \eLpx{\sum_{I\in\dyad{}} x_I h_I}^p\leq C^p\sum_{n=1}^{\infty}\eLpx{\sum_{I\in\dyad{}_{n-1}}  x_I h_I}^p.
\]
Since $\dyad{}_{n-1}$ consists of pairwise disjoint dyadic intervals we have 
\[
 \eLpx{\sum_{I\in\dyad{}_{n-1}}  x_I h_I}^p =\sum_{I\in\dyad{}_{n-1}}\nx{x_I}^p|I|.
\]
Consequently the inequality \eqref{eq:dyadtype} holds true if and only if 
\begin{equation}\label{def: haartype}
\eLpx{\sum_{I\in\dyad{}} x_I h_I}\leq C \left(\sum_{I\in\dyad{}}\nx{x_I}^p|I|\right)^{\frac{1}{p}}.  
\end{equation}

 Pietsch and Wenzel in their book \cite{Pietsch} on orthonormal systems and Banach space geometry introduced the notion of Haar type $p$ of a Banach space $X$ (see specifically \cite[Chapter 7]{Pietsch}). Recall that a Banach space $X$ is said to satisfy Haar type $p$ if there exists $C<\infty$ such that for any sequence $\{x_I\}_{I\in\mathcal{I}}\subseteq X$ (such that $\sum_{I\in\mathcal{I}} x_I h_I$ converges in $L^p(X)$)  the inequality \eqref{def: haartype} holds true. In that case we say that $X$ satisfies $HT_p$ and by $HT_p(X)$ we denote smallest constant $C$ such that \eqref{def: haartype} is satisfied. Thus Pisier's theorem \cite[Proposition 2.4 b) and Theorem 3.1 b)] {Pisier1975}(see also \cite[Theorem 10.22]{Pisier2016}) asserts that a Banach space X is of martingale type $p$ iff it is of Haar type $p$.

 Recall that we say $(\Omega,\, \FF,\, \left(\FF_n\right),\, \Pp)$ is a filtered probability space if $\left(\FF_n\right)$ is an increasing sequence of purely atomic sub-$\sigma$-algebras in a probability space $(\Omega,\,\FF,\,\Pp)$ and $\FF=\sigma\left(\bigcup \FF_n\right)$.
In the present paper we identify precisely all filtered probability spaces $(\Omega,\, \FF,\, \left(\FF_n\right),\, \Pp)$ that are able to determine the martingale type of Banach space $X$. We associate explicit intrinsic conditions on the filtration $(\FF_n)$ which determine that the upper $\ell^p$ estimates \eqref{type} imply the martingale type $p$ of Banach space $X$. As a consequence we obtain the following dichotomy.
\begin{T1}\label{main}
Let $1<p\leq 2$. For each fixed $(\Omega,\, \FF,\, \left(\FF_n\right),\, \Pp)$ the following dichotomy holds true:
Either, there exists $C>0$ such that for any Banach space $X$ and any $f\in L^p(\Omega,\, \FF,\, \Pp,\,X)$
\begin{equation}\label{filteredtype}
\eLpx{f}\leq C\left(\chng{\eLpp{\EE(f|\FF_0)}{\Omega,\, X}{p}^p}+ \sum_{n=1}^{\infty}\eLpx{\mardif{f}{n}}^p   \right)^{\frac{1}{p}},
\end{equation}
or the filtered probability space $(\Omega,\, \FF,\, \left(\FF_n\right),\, \Pp)$ and the upper $\ell^p$ estimates \eqref{filteredtype} already determine that the Banach space $X$ is of martingale type $p$.
\end{T1}
Theorem~\ref{infincarl}  of Section \ref{martineqz} and Theorem~\ref{thmfincarl} of Section \ref{sec:fincarl} are the main results of this paper. The dichotomy formulated in Theorem \ref{main} is a direct consequence thereof.

Next we point out the connection between Haar type and Carleson constants. In \cite{Geiss2008} the authors fixed a collection of dyadic intervals $\cgc\subset \mathcal{I}$ and obtained intrinsic conditions on the size of $\cgc$ such that  the upper $\ell^p$ estimates
\begin{equation}\label{dyadictypecol}
\eLpx{\sum_{I\in\cgc} x_I h_I}^p\leq C\sum_{I\in\cgc}\nx{x_I}^p|I|,
\end{equation}
for $x_I\in X$, $I\in\cgc$, imply that the Banach space $X$ is of Haar type $p$. Specifically in \cite{Geiss2008} it was shown that $\cgc$ and the upper $\ell^p$ estimates \eqref{dyadictypecol} determine the Haar type of $X$ if and only if the Carleson constant of $\cgc$ is unbounded, that is, 
\[
\crl{\cgc}:=\sup_{I\in\cgc} \sum_{\substack{J\in\cgc\\J\subseteq I }} \frac{|J|}{|I|}=\infty.
\]
The present paper could be viewed as continuation of \cite{Geiss2008}.

In \cite{muellerstud} the size of collections $\mathcal{G}\subset\dyad{}$, expressed in term of its Carleson constant $\crl{\mathcal{G}}$ was utilized in the solution of the Gamlen-Gaudet problem for dyadic $H^1$.

\section{Notation, conventions and definitions}\label{sec: not}
 Let $(\Omega,\, \FF,\, \left(\FF_n\right),\, \Pp)$ be a filtered probability space, where $\FF_n$, $n\geq 0$, is a $\sigma$-algebra generated by a finite family of atoms $\AAA_{n}$ satisfying $\Omega=\bigcup \limits_{A\in\AAA_n}A$ and $\AAA_0=\{\Omega\}$. Without changing the $p$-variation function
 \[
 \omega \rightarrow \nx{\EE( f|\FF_0)}^p+ \sum_{n=1}^{\infty} \nx{\mardifx{f}{n}{\omega}}^p
 \]
whose expectation appears in the statement of Theorem \ref{main} we may modify the
fields $\FF_n$ in such a way that every atom $A \in \AAA_n$ will be, either split up in strictly smaller atoms in $\AAA_{n+1}$, or else never split at any later step $m > n$. \chng{This change in the fields $\FF_n$ will clearly not modify the family $\bigcup_n A_n$ of all atoms, nor the union $\bigcup_n\FF_n$; also, it will not affect the definition
given below for the family $\Ec$ and the value of its Carleson constant $\crl{\Ec}$.}

\begin{figure}[h]
    \centering
    \includegraphics[width=\textwidth]{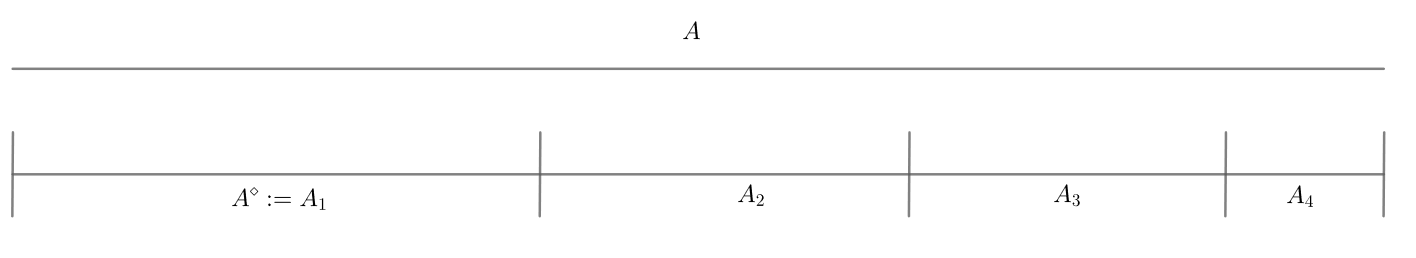}
    \caption{Here we put $A=A_1\cup A_2 \cup A_3\cup A_4$ and $N(A)=4$. Atoms are lined up according to their size from left to right.} In all our figures we adhere to the convention illustrated above.  
    \label{fig:int}
\end{figure}

Let
\[
\AAA=\bigcup_{n=0}^{\infty} \AAA_{n}.
\]
For $A\in \AAA_{n}$ we have $A=\bigcup\limits_{j=1}^{\Nay{A}} A_j$, where $A_j\in\AAA_{n+1}$ and are ordered in such a way that
\begin{equation}\label{notorder}
 \Pp(A_j)\geq \Pp(A_{j+1})>0     
\end{equation}

for every $j\in\{1,\ldots, \Nay{A}-1\}$ (see Figure \ref{fig:int}). For every $A\in\AAA$ we put $\st{A}:=A_1$.

Thus we arrive at the following critical definitions. We have $\AAA_0=\{\Omega\}.$ We put $\ccc_0=\emptyset$ and $\Ec_0=\AAA_{0}\backslash\ccc_{0}=\AAA_0$. We define for $n\in\NN\cup\{0\}\footnote{In this article we put $\NN=\{1,2,3,...\}$}$
\begin{equation*}
    \ccc_{n+1}=\{\st{A}: A\in \AAA_n\},\qquad\Ec_{n+1}=\AAA_{n+1}\backslash\ccc_{n+1}
\end{equation*}
and
\begin{equation}\label{defcole}
\ccc=\bigcup_{n=1}^{\infty} \ccc_{n}, \qquad\Ec=\bigcup_{n=0}^{\infty} \Ec_{n}.
\end{equation}
Note that $\Ec=\bigcup_{A\in\AAA}\{A_j: j\in\{2,\ldots,\Nay{A}\}\}\cup \{\Omega\}$.

Let $\Rc\subseteq \bigcup_{n} \FF_n$ be a nested collection, that is for $A$, $B\in\Rc$ such that $B\cap A\neq\emptyset$ we have $A\subseteq B$ or $B \subseteq A$. Given $A\in\FF$ we introduce the following subfamily
\[
\GC{1}{\Rc}{A}= \{ R : R\in \Rc,\; R\subsetneq A,\; R\mbox{ is maximal element with respect to inclusion}\},  
\]
and inductively we define for $k\in\{2,3, \ldots\}:$
\[
\GC{k}{\Rc}{A}=\bigcup\limits_{J\in \GC{k-1}{\Rc}{A}} G_{1}(\Rc,J). 
\]
Notice that collection $\Ec$ satisfies the following property: If $B\in \GC{k}{\Ec}{A}$ we have
\begin{equation}\label{eq: geodecrease}
   \Pp(B)\leq 2^{-k}\Pp(A). 
\end{equation}

For any collection of measurable sets $\Rc$ we define the point set covered by $\Rc$:
\[
\Est{\Rc}= \bigcup_{A\in \Rc} A.
\]
We will denote by $\crl{\Rc}$ the Carleson constant  of the collection $\Rc$ i.e.
\begin{equation}\label{defcrl}
\crl{\Rc}= \sup_{I\in \Rc}\frac{1}{\Pp(I)}\sum_{\substack{J\subseteq I\\ J\in \Rc}}\Pp(J).
\end{equation}

In \cite{Muellerisr} the collection $\Ec$ and its size, expressed in term of its Carleson constant $\crl{\Ec}$, was utilized in the classification theorem for martingale $H^1$ spaces.

\section{Infinite Carleson constant and Haar type}\label{martineqz}
 Recall that to a filtered probability space $(\Omega,\, \FF,\, \left(\FF_n\right),\, \Pp)$ we associated the collection $\Ec$ in equation \eqref{defcole} and its Carleson constant $\crl{\Ec}$ in equation \eqref{defcrl}. 
\begin{T1}\label{infincarl}
Let $(\Omega,\, \FF,\, \left(\FF_n\right),\, \Pp)$ be a filtered probability space, where each $\FF_n$ is purely atomic. Let $1<p\leq 2$ and $\crl{\Ec}= \infty$. Let $X$ be a Banach space. If there exists $T_p>0$ such that for any $f\in L^p(\Omega,X)$
\begin{equation}\label{typeX}
\eLpx{f}\leq T_p\left(\chng{\eLpp{\EE(f|\FF_0)}{\Omega,\, X}{p}^p}+ \sum_{n=1}^{\infty}\eLpx{\mardif{f}{n}}^p   \right)^{\frac{1}{p}},   
\end{equation}
 then the space $X$ is of Haar type $p$. Moreover 
 \begin{equation}\label{stala}
  HT_p(X)\leq 4^{\frac{2}{p}-1} T_p,     
 \end{equation}
 where $HT_p(X)$ is the Haar type p constant of $X$.
\end{T1}
\begin{Rem}
The constant on the right hand side of \eqref{stala} reflects the constant in Lemma \ref{Hfunc} H7). This refers also to the factor on the right side in \eqref{eq:varbou}.
\end{Rem}
Under the hypothesis $\crl{\Ec}=\infty$ we will show that $\operatorname{span}\{\chs_{A}\}_{A\in\Ec}$ supports/contains systems of functions equivalent to the initial segments of the Haar system in the Bochner-Lebesgue space $L^p(\Omega,X)$. Our approach may be divided into three separate steps as follows. First  in \mbox{Lemma \ref{Hfunc}} we design a "greedy algorithm" to  
construct a family of proto-Haar functions $\khan{A}$ associated to the collections of pairwise disjoint atoms $\GC{k}{\Ec}{A}$. Our greedy algorithm ensures that the following variational estimates hold true
\begin{equation}\label{eq:varbou}
\int_{\Omega}\sum_{n=1}^{\infty}|\mardifx{\khan{A}}{n}{\omega}| \dPx{\omega}\leq 4\Pp(A).    
\end{equation}
Definition \ref{ecsupph} and Lemma \ref{AAAimpMar} form the second part of the proof, where we explain the significance of the Gamlen-Gaudet construction in determining Haar type of a Banach space. Finally in Lemma \ref{unicolAAA} we show that $\crl{\Ec}=\infty$
 implies that $\Ec$ supports all initial segments of a generalized Haar system.
 
 We organize the proof in such a way that it is a direct consequence of Lemma~\ref{AAAimpMar} and Lemma~\ref{unicolAAA}.

\subsection{The Greedy Algorithm determines the selection of signs}
In this section we design the greedy algorithm that forms the main novelty of the present paper. It is employed in the proof of Lemma \ref{Hfunc} were we construct 
those proto-Haar functions $\khan{A}$ for which the variational estimates \eqref{eq:varbou} hold true. 
\begin{L1}\label{Hfunc}
 Let $0<\varepsilon<\frac{1}{2}$ and $k> -\ln_2(\varepsilon)$. Let $A\in \AAA$. If
\begin{equation}\label{larcon}
  \Pp(\Est{\GC{k}{\Ec}{A}})>\left(1-\frac{\varepsilon}{2}\right)\Pp(A)  
\end{equation}
holds true then there exists a function $\khan{A}$ such that 
\begin{enumerate}[label=H\arabic*), font=\normalfont]
    \item $|\khan{A}|\leq 1 \qquad \mbox{a.e.}$
    \item $\supp{\khan{A}}\subseteq A$,
    \item there exists $t:=t(A,\varepsilon, k, (\FF_n))\in\NN$ such that $\khan{A}$ is $\FF_t$ measurable.
    \item $\int_{A} \khan{A}\dP = 0$,
    \item $\Pp(\{|\khan{A}|\neq 1\}\cap A)\leq \varepsilon \Pp(A)$,
    \item $\left(\frac{1}{2}-\varepsilon\right)\Pp(A)\leq\Pp(\{\khan{A}=\xi\})\leq \frac{1}{2}\Pp(A)\qquad \forall\,\xi\in \{-1,1\}$,
    \item $\sum\limits_{n=1}^{\infty}\eLpp{\mardif{\khan{A}}{n}}{\Omega}{1} \leq 4 \Pp(A)$.
\end{enumerate}
\end{L1}

\begin{proof}
   Let $m\in\NN\cup\{0\}$ and $A\in \AAA_m$. We will construct sequences of sets $(Y_l)_{l=m}^{\infty}$, $(U_l)_{l=m}^{\infty}$, $(V_l)_{l=m}^{\infty}$ and functions $(\khan{l})_{l=m}^{\infty}$ such that for every $l\geq m$ we have
\begin{enumerate}
    \item $\khan{l}$ satisfies properties H1), H2), H4), H7),
    \item $\{|\khan{l}|\neq 1\}\cap A= Y_l\in\AAA_l$,
    \item $\{\khan{l}= 1\}\cap A= U_l\in\FF_l$,
    \item $\{\khan{l}= -1\}\cap A= V_l\in\FF_l$,
    \item $\Pp(U_l) < \frac{1}{2}\Pp(A)\leq  \Pp( Y_{l})+ \Pp(U_{l})$
    \item $Y_{l+1}\subset Y_{l}$,
    \item $U_l\subset U_{l+1}$,
    \item $V_l\subset V_{l+1}$,
    \item $U_{l} \cup V_{l} \cup Y_{l}= A$,
    \item $H_l$ is $\FF_l$-measurable.
\end{enumerate}

    We put $Y_m=A$, $H_m\equiv 0$, $U_m=V_m=\emptyset$ and $w_A=\frac{1}{2}\Pp(A)$. List of properties is obviously satisfied.
  Assume that for $l\geq m$ we have constructed $Y_l$, $U_{l}$, $V_l$, $H_l$ satisfying the properties listed above. Recall that $Y_l\in\AAA_l$. Hence  $Y_l=\bigcup_{j=1}^{N(Y_{l})} Y_{l,j}$, where $Y_{l,j}\in\AAA_{l+1}$ are disjoint atoms and are ordered according to their size (see \eqref{notorder}). Since
    \[
    \Pp(U_l) < w_A\leq  \Pp( Y_{l})+ \Pp(U_{l}),
    \]
    there exists an index $s$ such that
     \[
    \Pp(U_l)+\sum_{j<s} \Pp(Y_{l,j}) < w_A\leq  \Pp( Y_{l,s})+ \Pp(U_{l}) +\sum_{j<s} \Pp(Y_{l,j}),
    \]
    We define 
    \[
        U_{l+1}=U_l \cup \bigcup_{j<s} Y_{l,j},\qquad V_{l+1}=V_l \cup \bigcup_{j>s} Y_{l,j}\quad\mbox{and}\quad Y_{l+1}=Y_{l,s}.
    \]
    Observe that $U_{l+1} \cup V_{l+1} \cup Y_{l+1} = U_{l} \cup V_{l} \cup Y_{l}= A$. Moreover
    \[
     \Pp(U_{l+1}) - \Pp( Y_{l+1}) - \Pp(V_{l+1})= 2 \Pp(U_{l+1}) -\Pp(A) < 0.
    \]
    and
    \[
    0 \leq 2 \Pp(U_{l+1}) + 2 \Pp( Y_{l+1})-\Pp(A)= \Pp(U_{l+1}) + \Pp( Y_{l+1}) - \Pp(V_{l+1})
    \]
    Therefore there exists $-1<c_{l+1}\leq 1 $ such that
    \[
    \Pp(U_{l+1}) + c_{l+1} \Pp( Y_{l+1}) - \Pp(V_{l+1})=0.
    \]
    We define 
    \[
    \khan{l+1}= \chs_{U_{l+1}}+ c_{l+1} \chs_{Y_{l+1}}-\chs_{V_{l+1}}.
    \]
    Properties of sets $Y_l$, $U_{l}$, $V_l$ are obvious by the construction. Function $\khan{l+1}$ clearly is $\FF_{l+1}$-measurable. In fact properties H1), H2), H4) of $\khan{l+1}$ are obvious from the definition of $\khan{l+1}$. Only property H7) demands more attention.

    Observe that the sequence of functions $(\khan{n})$ forms a martingale with respect to the filtration $(\FF_n)$. Indeed functions $\khan{l}$ and $\khan{l+1}$ are equal outside of the set $Y_l\in \AAA_l$ and
    \[
    \EE \khan{l} \chs_{Y_l}= - \EE \khan{l} \chs_{\Omega\backslash Y_l}= - \EE \khan{l+1} \chs_{\Omega\backslash Y_l}=\EE \khan{l+1}  \chs_{Y_l}.
    \]
    We calculate that
    \[
    0=\EE(\khan{l+1}-\khan{l})= (1-c_l) \Pp(U_{l+1}\backslash U_{l}) + (c_{l+1}-c_{l}) \Pp(Y_{l+1}) + (-1-c_{l}) \Pp(V_{l+1}\backslash V_{l}).
    \]
    Therefore
    \begin{equation}\label{eq: 11206}
    \begin{split}
       |(c_{l}-c_{l+1}) \Pp(Y_{l+1})|
       &\leq |(1-c_l) \Pp(U_{l+1}\backslash U_{l})|+|(-1-c_{l}) \Pp(V_{l+1}\backslash V_{l})|
       \\&\leq 2 \left( \Pp(V_{l+1}\backslash V_{l})+\Pp(U_{l+1}\backslash U_{l})\right).      
    \end{split}
    \end{equation}
    Observe that
    \[
    \sum\limits_{n=0}^{\infty}\eLpp{\EE(\khan{l}\big|\FF_{n+1})-\EE(\khan{l}\big|\FF_{n})}{\Omega}{1}= \sum\limits_{n=m}^{l-1}\eLpp{H_{n+1}-H_{n}}{\Omega}{1}.
    \]
    Recall that $(U_l)$, $(V_l)$ are sequences of increasing sets. It follows from this and \eqref{eq: 11206} that 
    \[
    \begin{split}
        \sum\limits_{n=m}^{l-1}\eLpp{H_{n+1}-H_{n}}{\Omega}{1}
         &= \sum\limits_{n=m}^{l-1}|1-c_n| \Pp(U_{n+1}\backslash U_{n})+|c_{n}-c_{n+1}| \Pp(Y_{n+1})
         +|1+c_{n}| \Pp(V_{n+1}\backslash V_{n})
         \\&\leq 4 \sum\limits_{n=m}^{l-1} \left( \Pp(V_{n+1}\backslash V_{n})+\Pp(U_{n+1}\backslash U_{n})\right)
         \\&\leq  4 \left( \Pp(V_{l})+\Pp(U_{l})\right)
         \\&\leq  4 \Pp(A).
    \end{split}
    \]
    
Now we will prove the existence of an index $t$ such that $H_t$ also satisfies H5) and  H6). Then we will put $H_A:=H_t$. By \eqref{larcon} there exists a finite family of disjoint atoms $\{B_j\}_{j=1}^{s} \subset \GC{k}{\Ec}{A}$ such that
\begin{equation}\label{eq: 146}
   \sum_{j=1}^{s}\Pp(B_j)\geq (1-\varepsilon) \Pp(A) 
\end{equation}
Since the family is finite we can find $t$ such that $\{B_j\}_{j=1}^{s}\subset\FF_{t-1}$. Since $Y_t\in\AAA_t$ and $\{B_j\}\subset\AAA\cap\FF_{t-1}$ either there exists $j$ such that $Y_t\subset B_j$ or $Y_t$ and $B_j$ are disjoint for every $j$. If  $Y_t\subset B_j$ then
\[
\Pp( \{|\khan{t}|\neq 1\}\cap A)= \Pp(Y_{t})\leq \Pp(B_j) \stackrel{\eqref{eq: geodecrease}}{\leq} 2^{-k} \Pp(A) \stackrel{k>-\ln\varepsilon}{\leq} \varepsilon \Pp(A).
\]
Otherwise $Y_t$ is disjoint from the union of sets $B_j$ and by \eqref{eq: 146}
\[
\Pp( \{|\khan{t}|\neq 1\}\cap A)\leq \Pp(A)- \sum_{j=1}^{s}\Pp(B_j)\leq \varepsilon\Pp(A).
\]
This proves that $H_t$ satisfies H5). Now for H6) we only need to observe that from the construction the following estimate is satisfied: 
\[
w_A -\Pp(Y_t)\leq \min\{\Pp(U_{t}),\Pp(V_{t})\}\leq \max\{\Pp(U_{t}),\Pp(V_{t})\}\leq w_A.
\]

\end{proof}

\begin{Rem}\label{MREX}
The proto-Haar function $\khan{A}
$ was constructed by means of the greedy algorithm presented in the proof of Lemma \ref{Hfunc}. As a critical consequence the variational norm \eqref{eq:varbou} of the function $\khan{A}$ is bounded. See Lemma \ref{Hfunc} condition H7).
We point out that the complexity of our construction can not be reduced entirely. This became especially clear in view of the following example communicated to us by Maciej Rzeszut. We are grateful for his permission to present it here. We fix  $0<\varepsilon<\frac{1}{2}$. There \chng{exists} a function $\tilde{H}_{[0,1]}$ such that for the dyadic filtration it satisfies H1)-H6) but it fails to satisfy H7). We select $N_0\in\NN$ such that for any $N>N_0$ 
\[
\Pp\bigg(\{\omega : \Theta_{N}(\omega)=0\}\bigg)\leq \varepsilon,
\]
where 
\[ \Theta_N(\omega)=\sum_{I\in \dyad{N}} h_I(\omega),
\]
$\dyad{N}=\{I\;:\;I\mbox{ is dyadic and }|I|\geq 2^{-N} \}$ and $\{h_I\}_{I\in\dyad{N}}$ is the initial segment of Haar system. Consider the $\{0,1,-1\}$ valued function 
\begin{equation}\label{201001}
f(\omega)=\operatorname{sgn}\left(\Theta_N(\omega)\right),
\end{equation}
We shall employ Khintchine's inequality to see that
\begin{equation}\label{201002}
\int_{\Omega}\sum_{n=1}^{N+1}|\mardifdyad{f}{n}| \dP \geq \sqrt{\frac{N+1}{2}},    
\end{equation}
where $\dyadsigma{n}$ denotes the $\sigma$-algebra generated by $\dyad{n}$. We write $\mardifold{g}{n}=\mardifdyad{g}{n}$ for a given function $g$. Recall also $\mathcal{I}_{n}=\dyad{n}\backslash\dyad{n-1}$, $\mathcal{I}_0=\{[0,1]\}$ and that $r_n=\sum_{I\in\mathcal{I}_{n}} h_I$ are independent Rademacher variables.
\[
\begin{split}
 \sqrt{\frac{N+1}{2}}&\stackrel{Khintchine}{\leq}\EE\left|\sum_{n=0}^{N} r_n\right|
= \EE\left(\operatorname{sgn}\left(\sum_{n=0}^{N} r_n\right)\sum_{n=0}^{N} r_n\right)
\\&=\EE\left(f\sum_{n=0}^{N} r_n\right)
=\sum_{n=0}^{N}\EE\left( f \mardifold{r_n}{\chng{n+1}}\right)
\\&=\sum_{n=0}^{N}\EE\left( r_n\mardifold{f}{\chng{n+1}}\right)
\leq\sum_{n=0}^{N}\EE\left( \left|r_n \right|\left|\mardifold{f}{\chng{n+1}}\right| \right)
\\&=\sum_{n=0}^{N} \EE\left( \left|\mardifold{f}{\chng{n+1}}\right|\right)
=\sum_{n=0}^{N} \eLpp{\mardifold{f}{\chng{n+1}}}{\Omega}{1}\chng{.}
\end{split}
\]
In summary each of the functions $f_N=\operatorname{sgn}(\Theta_N)$ satisfies H1)-H6) and for any constant $C$ in H7) condition H7) is violated for $N$ large enough.
\end{Rem}
\begin{Rem}
The extremal nature of the lower bound \eqref{201002} follows from H\"{o}lder's inequality and Parseval's equation. Indeed for any $g\in L^{2}(\Omega)$ satisfying $\eLpp{g}{\Omega}{2}=1$ we have
\[
\begin{split}
\int_{\Omega}\sum_{n=1}^{N}|\mardifdyad{g}{n}| \dP&\leq\int_{\Omega}\sqrt{N}\left(\sum_{n=1}^{N}|\mardifdyad{g}{n}|^2\right)^{\frac{1}{2}}\dP
\\&\leq \sqrt{N} \left(\int_{\Omega}\sum_{n=1}^{N}|\mardifdyad{g}{n}|^2\dP\right)^{\frac{1}{2}}
\\&\leq\sqrt{N}\eLpp{g}{\Omega}{2}=\sqrt{N}.
\end{split}
\]
\end{Rem}

\subsection{Gamlen and Gaudet determine Haar type}
In this subsection for any $n$, similarly as in \cite{Geiss2008}, we introduce a family of conditions, which will imply Haar type. Recall that $\dyad{n}=\{I\;:\;I\mbox{ is dyadic and }|I|\geq 2^{-n} \}$.
\begin{Def}\label{ecsupph}
For a filtered probability space $(\Omega,\, \FF,\, \left(\FF_n\right),\, \Pp)$ we say that the corresponding collection $\Ec$ supports all initial segments of a generalized Haar system if there exists a constant $K>0$ such that for any $n\in\NN$, $\delta>0$ and $\varepsilon=\delta 2^{-n-1} $ there are a family of collections $\{\col{I}\}_{I\in\dyad{n}}$, a family of functions $\{g_I\}_{I\in\dyad{n-1}}$ and an atom  $\Io\in\Ec$ satisfying the following conditions:  
\begin{enumerate}[label= $\AAA$\!\arabic*)]
    \item $\col{[0,1]}=\{\Io\}.$
    \item $\col{I}\subsetneq\{A\in\Ec\;: A\subset \Io\}$ for $I\in\dyad{n}.$
    \item The elements of $\col{I}$ are pairwise disjoint for every $I\in \dyad{n}$.
 
    \item $\Est{\col{I}}\subseteq \Est{\col{J}}$ if and only if $I\subseteq J$ and
    $\Est{\col{I}}\cap \Est{\col{J}}\neq \emptyset$ if and only if $I\cap J\neq \emptyset$.
     \item For any $I\in\dyad{n-1}$ we have  $g_I:\Omega\rightarrow [-1,1]$,
    \[
    \supp{g_I}\subseteq \Est{\col{I}}
    \]
    and $\Est{\col{I^+}} \subseteq \{g_I=1\}$, $\Est{\col{I^-}} \subseteq \{g_I=-1\}$, where $I^-$ is the right half of $I$ and $I^{+}$ is the left half of $I$. 
       \item There exists $S\in \NN$ such that the elements of $\col{I}$ are $\FF_{S}$ measurable for every $I\in\dyad{n}$.
    
    \item Functions $g_I$ are $\FF_{S}$ measurable for any $I\in\dyad{n-1}$.
    \item For any atom $A\in \AAA_{S}$ and any $m\in\{1,\ldots, S\}$ there exists at most one dyadic interval $I\in\dyad{n-1}$ such that
    \[
    \left( \mardif{g_I}{m}\right)\big|_A\neq 0.
    \]
    We will denote that interval by $I(A,m)$ (Remark \ref{adnotationA8} below restates a combinatorially dual picture of this condition).
    \item For any $I\in\dyad{n-1}$ 
    \begin{equation}\label{margijaki1}
  \eLpp{\EE(g_I|\FF_0)}{\Omega}{1} +  \sum_{j=1}^{\infty} \eLpp{\mardif{g_I}{j}}{\Omega}{1} \leq K\Pp(\Est{\col{I}}).
    \end{equation}
    \item For every $I\in \dyad{n}$ we have
    \begin{equation}\label{likedyad}
    (|I|-2\varepsilon)\Pp(\Io)\leq \Pp(\Est{\col{I}})\leq |I|\Pp(\Io) .
    \end{equation}
    
\end{enumerate}

\end{Def}
\begin{Rem}
    By $\AAA$\!10) for $|I|=2^{-n}$ we have 
    \[
    \frac{1-\delta}{2^{n}}\Pp(\Io)\leq \Pp(\Est{\col{I}})\leq \frac{1}{2^{n}}\Pp(\Io) \]
    and
    \[
    (1-\delta)\Pp(\Io)\leq  \sum_{|I|=2^{-n}} \Pp(\Est{\col{I}})\leq \Pp(\Io).
    \]
\end{Rem}
\begin{Rem}
Note that the condition $\AAA 9)$ implies that for $1<p<2$ and $\theta_p=2-p$ we have
  \begin{equation}\label{margijakip}
   \sum_{j=1}^{\infty} \eLp{\mardif{g_I}{j}}^p \leq K^{\theta_p}\Pp(\Est{\col{I}}).
    \end{equation}
Indeed H\"{o}lder's inequality yields 
\[
\begin{split}
\sum_{j=1}^{\infty}\chng{\eLpp{\mardifold{g_{I}}{j}}{\Omega}{p}^p}&\leq\left( \sum_{j=1}^{\infty}\eLpp{\mardifold{g_I}{j}}{\Omega}{1}\right)^{\theta_p}\left(\sum_{j=1}^{\infty}\eLpp{\mardifold{g_I}{j}}{\Omega}{2}^2\right)^{1-\theta_p}
\\&\stackrel{\eqref{margijaki1}}{\leq} K^{\theta_p}\Pp(\Est{\col{I}}),
\end{split}
\]
where $\mardifold{f}{j}=\mardif{f}{j}$.
\end{Rem}
\begin{Rem}\label{adnotationA8}
In fact our functions $\{g_I\}$ satisfy the following conditions. Given an atom $A\in\AAA_S$ and given two dyadic intervals $I\subsetneq J$ there exist disjoint integer intervals $N_I\subseteq\{1,\ldots, S\}$ and $N_J\subseteq\{1,\ldots, S\}$ such that
\[
\sup  N_J < \inf N_I
\]
and
\[
g_I|_A = \sum_{j\in N_I} \mardifold{g_I}{j}|_A, \qquad\mbox{ and }\qquad g_J|_A = \sum_{j\in N_J} \mardifold{g_J}{j}|_A.
\]
\end{Rem}

Now we turn to showing that conditions $\AAA1)-\AAA10)$ imply Haar type.
\begin{L1}\label{AAAimpMar}
Let $1< p\leq 2$. Assume that the collection $\Ec$ supports all initial segments of a generalized Haar system. Then the Banach space $X$ satisfies inequality \eqref{typeX} iff $X$ has Haar type $p$.
\end{L1}
\begin{proof}
 Let $n\in \NN$ and $0<\delta<1$. We fix a finite sequence $(x_J)_{J\in\dyad{n-1}}$ of vectors from Banach space $X$. In view of Definition \ref{ecsupph} there are collections $\col{J}\subsetneq \Ec$ with $J\in\dyad{n}$, functions $g_J$ with $J\in\dyad{n-1}$ and $\Io\in\Ec$ satisfying \mbox{$\AAA$1)-$\AAA$10)}.

For fixed dyadic interval $I$ of length $2^{-n}$ we choose $U_I\subset V_I\subset I$ with $|U_I|= (|I|-2\varepsilon)$ and $|V_I|=\frac{\Pp(\Est{\col{I}})}{\Pp(\Io)}.$ Next define
\[
\Wb=\bigcup_{|I|=2^{-n}} \Est{\col{I}},\qquad \Wv=\bigcup_{|I|=2^{-n}} V_{I},\qquad \Wu=\bigcup_{|I|=2^{-n}} U_{I}.
\]
Let 
\[
f(\omega)= \sum_{J\in\dyad{n-1}} x_J g_J(\omega). 
\]
The function $f$ and each of the functions $g_J$, where $J\in\dyad{n-1}$, are constant on any of the sets $\Est{\col{I}}$ with $|I|=2^{-n}$. Indeed, either 
\[
\Est{\col{I}}\cap\Est{\col{J}}=\emptyset\quad\mbox{ or }\quad \Est{\col{I}}\subseteq \Est{\col{J^+}}\subseteq \{g_J=1\}\quad\mbox{ or }\quad \Est{\col{I}}\subseteq \Est{\col{J^-}}\subseteq \{g_J=-1\}.
\]
We define $\tilde{f}:\Wv\rightarrow X$ and $\tilde{g}_J :\Wv\rightarrow X$ :
\[
\tilde{g}_J(V_I)= g_J(\Est{\col{I}})\qquad \mbox{and}\qquad \tilde{f}= \sum_{J\in\dyad{n-1}} x_J \tilde{g}_J.
\]
Clearly for $|I|=2^{-n}$
\[
\tilde{f}(V_I)=f(\Est{\col{I}}).
\]
Observe that on the probability space $(\Wu, (1-\delta)^{-1} dt)$ the system $(\tilde{g}_J)_{J\in\dyad{n-1}}$ has the same joint distribution as the usual Haar basis $(h_J)_{J\in\dyad{n-1}}$ on the unit interval. Therefore
\[
\begin{split}
\left\|\sum_{I\in\dyad{n-1}} x_I h_I\right\|^p_{L^p([0,1],X)}&=(1-\delta)^{-1} \int_{\Wu} \nxb{\tilde{f}(x)}^p \dx \leq (1-\delta)^{-1} \int_{\Wv} \nxb{\tilde{f}(x)}^p \chng{\dx}.
\end{split}
\]
On the other hand one can easily express the integral of $\tilde{f}$  in terms of the function $f$. 
\[
\int_{\Wv} \nxb{\tilde{f}(x)}^p \dx=\frac{|\Wv|}{\Pp(\Wb)}\int_{\Wb} \nx{f}^p\dP\leq\frac{1}{\Pp(\Io)}\int_{\Omega} \nx{f}^p\dP.
\]

Given $j\in \{1,2,\ldots, S\}$ we write 
\[
\mathscr{X}_j=\{(A,J): A\in \AAA_S \mbox{ and } J=I(A,j)\}
\]
and for fixed $J\in\dyad{n-1}$ we put 
\[
\mathscr{Y}_J=\{(A,j): A\in \AAA_S \mbox{ and } J=I(A,j)\}.
\]
Recall that $I(A,j)$ is defined by $\AAA$8).  We use \eqref{typeX} to obtain
\[
\begin{split}
\int_{\Omega} \nx{f}^p\dP&\stackrel{\eqref{typeX}}{\leq}T_p^p\sum_{j=1}^{S}\eLpx{\mardif{f}{j}}^p
\\&\stackrel{\AAA 7)}{=}T_p^p\sum_{j=1}^S\sum_{A\in \AAA_S}\eLpxA{\mardif{f}{j}}{A}^p
\\&\stackrel{\AAA 8)}{=}T_p^p\sum_{j=1}^{S}\sum_{(A,J)\in\mathscr{X}_j}\eLpxA{\mardif{x_J g_{J}}{j}}{A}^p    
\\&=T_p^p\sum_{J\in \dyad{n-1}}\sum_{(A,j)\in\mathscr{Y}_J}\eLpxA{\mardif{x_J g_{J}}{j}}{A}^p 
\\&=T_p^p\sum_{J\in \dyad{n-1}}\sum_{j\in \NN}\nx{x_J}^p\eLp{\mardif{g_{J}}{j}}^p 
\\&\stackrel{\eqref{margijakip},\eqref{likedyad}}{\leq} T_p^p K^{\theta_p }\sum_{J\in \dyad{n-1}} \nx{x_J}^p |J|\Pp(\Io)
\\&= T_p^p K^{\theta_p}\sum_{J\in \dyad{n-1}} \nx{x_J}^p \eLpp{\haar{J}}{[0,1]}{p}^p\Pp(\Io).
\end{split}
\]
In summary we get 
\[
 \left\|\sum_{J\in\dyad{n-1}} x_J h_J\right\|^p_{L^p([0,1],X)}\leq \frac{T_p^p K^{\theta_p}}{(1-\delta)}\sum_{J\in \dyad{n-1}} \nx{x_J}^p \eLp{\haar{J}}^p.
\]
 The above estimate no longer depends on the choice of the collection $\{\col{I}\}$. Since $\delta$ was an arbitrary number from interval $(0,1)$ we let $\delta\rightarrow 0$:
\[
\left\|\sum_{I\in\dyad{n-1}} x_I h_I\right\|^p_{L^p([0,1],X)}\leq T_p^p K^{\theta_p}\sum_{J\in \dyad{n-1}} \nx{x_J}^p \eLp{\haar{J}}^p. 
\]
This means that $X$ has Haar type $p$ with constant
\begin{equation}\label{estmartconst}
HT_p(X)\leq T_p K^{\frac{\theta_p}{p}}=T_p K^{\frac{2}{p}-1}.    
\end{equation}
Recall that Pisier proved in \cite[Proposition 2.4 b) and Theorem 3.1 b)] {Pisier1975} that a Banach space of Haar type $p$ is also a Banach space of martingale type $p$. This means that for any filtered probability space $(\Omega,\, \FF,\, \left(\FF_n\right),\, \Pp)$ the inequality \eqref{type} is satisfied.

\end{proof}
\subsection{Carleson-Garnett condensation lemma}
Here we want to show that $\crl{\Ec}=\infty$ implies that $\Ec$ supports all initial segments of a generalized Haar system. In the proof of this fact Lemma \ref{Hfunc} plays a crucial role. In order to show that condition \eqref{larcon} from Lemma \ref{Hfunc} is satisfied we use that for collections with infinite Carleson constant the following  version of the condensation lemma holds true (see  \cite[Ch. X : Lemma 3.2 and p. 414]{Garnett}).
\begin{L1}\label{enhcond1}
Let $\Rc\subseteq\Ec$ be a nested collection such that $\crl{\Rc}=\infty$. For any $\tilde{\varepsilon}>0$, $n\in\NN$ and $k\in\NN$ there are families of sets $ \{\Tc_j\}_{j\in\{0,1,\ldots,n\}}$ and a set $\Io\in \Rc$ satisfying
\begin{enumerate}[label=(\roman*),font=\normalfont]
    \item $\Tc_0=\{\Io\}$,
    \item $\Tc_j\subseteq \GC{ k j}{\Rc}{\Io}$ for any $j\in\{1,\ldots,n\}$,
    \item For $j\in\{1,2,\ldots,n\}$ and $A\in\Tc_{j-1}$ we have
    \[
    \Pp(A\cap \Est{\Tc}_{j})>(1-\tilde{\varepsilon}) \Pp(A),
    \]
    \item For any $j\in\{1,2,\ldots,n\}$ the family $\Tc_j$ is finite.
\end{enumerate}
\end{L1}
\begin{proof}
Let us fix auxiliary notation. We define families of sets $(g_j(A))_{j\in\NN} $ by
\[
g_0(A)=\{A\},\quad\mbox{and}\quad g_l(A)=\GC{kl}{\Rc}{A}
\]
for $l>0$. Clearly $\Est{g_l(A)}\subset A$. The family $g_l(A)$ consists of pairwise disjoint sets. The sets $\Est{g_m(B)}$ are pairwise disjoint for $B\in g_l(A)$. Thus
\begin{equation}\label{244}
    g_{l+m}(A)= \bigcup_{B\in g_l(A)}g_m(B)\quad\mbox{and}\quad \Pp(\Est{g_{l+m}(A)})=\sum_{B\in g_l(A)}\Pp(\Est{g_m(B)})
\end{equation} 
Now we will prove an auxiliary lemma.
\
\begin{L1}\label{te328}
Let $\varepsilon\in (0, 1)$, let $p, m > 0$ be integers and suppose that
\[
\Pp(g_{m+1}(A))> (1-\varepsilon^{p+1})\Pp(A).
\]
Let 
\[
\Tc:=\Tc(A,m+1,p) = \{ B \in g_1(A) : \Pp(\Est{g_m(B)}) > (1 -\varepsilon^p)\Pp(B)\}.
\]
Then
\[
\Pp(\Est{\Tc}) >(1 -\varepsilon) \Pp(A).
\]   
\end{L1}
\begin{proof}
Let $a_0 := \Pp(A)$ and $a_1 := \Pp(g_{m+1}(A))\leq a_0$. By the assumptions, we have that
\begin{equation}\label{245}
    a_1>(1-\varepsilon^{p+1})a_0.
\end{equation}
Let
\[
\Ucc = g_1(A) \backslash \Tc = \{B \in g_1(A) : \Pp(\Est{g_m(B)})
\leq (1 -\varepsilon^p) \Pp(B)\}
\]
and
\[
t_0 := \Pp(\Est{\Tc}),\quad u_0 = \Pp(\Est{\Ucc}),\quad t_1 =\sum_{B\in\Tc}\Pp(\Est{g_m(B)})\leq t_0,\quad u_1 = \sum_{B\in \Ucc}\Pp(\Est{g_m(B)})\leq u_0.
\]
By the definition of $\Ucc$, we know that
\begin{equation}\label{255}
  u_1\leq\sum_{B\in\Ucc} (1-\varepsilon^p)\Pp(B)\leq (1-\varepsilon^p)\Pp(\Est{\Ucc})= (1-\varepsilon^p) u_0
\end{equation}
Using \eqref{245}, \eqref{244} and \eqref{255} we have
\begin{equation}\label{259}
 (1-\varepsilon^{p+1})a_0<a_1=\sum_{B\in g_1(A)}\Pp(\Est{g_m(B)})= t_1+u_1\leq t_1   +(1-\varepsilon^p) u_0 
\end{equation} 
Since $t_0 + u_0 = \Pp(\Est{g_1(A)})\leq \Pp(A)=a_0$ it follows that
\[
\varepsilon^p u_0 = (t_0 + u_0) - (t_0 + (1 -\varepsilon^p))u_0 < a_0 - (1 - \varepsilon^{p+1})a_0=\varepsilon^{p+1}a_0,
\]
namely $u_0 < \varepsilon a_0$, so that \eqref{259} implies 
\[ 
(1-\varepsilon^{p+1})a_0 < t_0 + (\varepsilon - \varepsilon^{p+1})a_0\quad\mbox{thus}\quad(1 -\varepsilon)a_0 < t_0.  
\]

\end{proof}

Now fix $\varepsilon=\frac{\tilde{\varepsilon}}{2}$. Since $\crl{\Rc}=\infty$, there exists a set $A$ such that:
\[
\sum_{\substack{B\subset A\\B\in\Rc}} \frac{\Pp(B)}{\Pp(A)}o > \frac{k n}{\varepsilon^{n+1}}.
\]
By \cite[Lemma 3.1.4]{Mueller2005} there exists $\Io\subset A$, $\Io\in\Rc$ with
\[
\sum_{B\in \GC{kn}{\Rc}{\Io}} \Pp(B)> (1-\varepsilon^{n+1})\Pp(\Io).
\]
We put $\Tc_0=\{\Io\}$. Using Lemma \ref{te328} we define inductively 
\[
\Tc_{j}:=\bigcup_{A\in\Tc_{j-1}} \Tc(A,n-(j-1),n-(j-1)) \qquad j\in\{1,2,\dots,n\}.
\]
The definition of the families $\Tc_j$ allows one to apply inductively Lemma \ref{te328}, the result of Lemma \ref{te328} at one step giving the assumption for the next step.   By the construction of $\Tc_j$ conditions $(i),(ii)$ are satisfied. For $A\in\Tc_{j-1}$ by Lemma \ref{te328} we get
\[
\Pp(A\cap \Est{\Tc}_{j})> (1-\varepsilon)\Pp(A).
\]
Therefore condition (iii) is satisfied for $\varepsilon=\frac{\tilde{\varepsilon}}{2}$. Since we have space for error we modify $\Tc_j$ to be finite families.
\end{proof}

We are ready for the next step towards the proof of Theorem \ref{infincarl}.
\begin{L1}\label{unicolAAA}
If $\crl{\Ec}=\infty$ then $\Ec$ supports all initial segments of a generalized Haar system with constant $K=4$.
\end{L1}
\begin{proof}
Fix $\frac{1}{2}>\varepsilon>0$ and $n\in\NN$. We choose $k$ and $\tilde{\varepsilon}$ such that $2^{-k}<\tilde{\varepsilon}<\frac{\varepsilon}{2}$. We apply the Lemma \ref{enhcond1} for $\tilde{\varepsilon}$, $k$ and $\Rc=\Ec$. We obtain a set $\Io\in \Ec$ and families $\Tc_j\subseteq \Ec$, where $j\in\{0,1,\ldots,n\} $  satisfying conditions (i)-(iv) of Lemma \ref{enhcond1}. We define $\col{[0,1]}=\{\Io\}$ and $g_{[0,1]}=\khan{\Io}.$ Now we define $\col{I}$ and $g_I$ inductively with respect to canonical ordering. Fix $I$ such that $|I|=2^{-l}$. Assume that $\col{I}$ and $g_I$ are already defined. Recall that $I^{+}$ denotes the left half of $I$ and $I^{-}$ the right one. We define
\[
\col{I^{+}}=\{B\in \Tc_{l+1}\;:\;B\subsetneq \{g_I=  1\} \},\qquad \col{I^{-}}=\{B\in \Tc_{l+1}\;:\;B\subsetneq \{g_I= - 1\} \}.
\] 
Note that $\Tc_{l+1}\subseteq\AAA$. Hence $\col{I^{+}},\col{I^{-}}\subseteq\AAA$ and by Lemma \ref{enhcond1} (ii),(iii) the assumptions of Lemma \ref{Hfunc} are satisfied for every atom $A\in \Tc_{l+1}$.  Thus for every $A\in \Tc_{l+1}$ the function $\khan{A}$ is well defined. We put
\[
g_{I^{+}}=\sum_{A\in\col{I^{+}}} \khan{A},\qquad g_{I^{-}}=\sum_{A\in\col{I^{-}}} \khan{A}.
\]
It is easy to observe that properties $\AAA1)-\AAA5)$ are satisfied. By Lemma \ref{enhcond1} iv) and Lemma \ref{Hfunc} H3) we know that properties $\AAA 6)$, $\AAA 7)$ hold.

Now we prove $\AAA 8)$.  We fix an atom $A\in \AAA_{S}$ and $m\in\{1,\ldots, S\}$. We show that for any pair of distinct dyadic intervals $I,\;J\in\dyad{n-1}$ we have
\[
    \left(\mardif{g_I}{m}\right)\big|_A= 0\qquad \mbox{or}\qquad
    \left(\mardif{g_J}{m}\right)\big|_A= 0.
 \]
By $\AAA4)$ and $\AAA5)$ it is enough to consider pairs of dyadic $I$, $J$ satisfying $I\subsetneq J$. It suffices to consider $A\subseteq \Est{\col{I}}$. Next determine $B_1$ such that $B_1\in \col{I}$ and $A\subseteq B_1$. Let
\[
\mathcal{V}=\{B\in\AAA\;:\;B_1\subseteq B\mbox{ and } g_J\mbox{ is constant on  }B\}.
\]
Note that $\Est{\col{I}}$ is contained in $\Est{\col{J^+}}$ or in $\Est{\col{J^-}}$, hence $B_1 \subseteq C \in \col{J^+}$ or $B_1 \subseteq C \in \col{J^-}$ for some atom $C\in\AAA$. Note that $\mathcal{V}$ is non empty set. Let $B_2$ be a maximal element of family $\mathcal{V}$ with respect to inclusion.
There exists a unique $l$ satisfying $B_2\in \AAA_l$. If $l < m$ we have
\[
\EE(g_J|\FF_m)|_A=\EE(g_J|\FF_{m-1})|_A.
\]
Therefore
\[
\left(\mardif{g_J}{m}\right)\big|_A=0.
\]
If conversely $l\geq m$ then
\[
\left(\mardif{g_I}{m}\right) |_A=0.
\]
Indeed by Lemma \ref{Hfunc} H4) we have for $l\leq m $
\[
\EE(\khan{B_1}|\FF_m)\big|_A=0\qquad\mbox{and}\qquad\EE(\khan{B_1}|\FF_{m-1})\big|_A=0.
\]
In summary for fixed $A\in\AAA_{S}$ and $m\in\{1,\ldots, S\}$ there is at most one dyadic interval $I$ satisfying
\[
\mardifold{g_I}{m}\big|_A=\mardif{g_I}{m}   \big|_A\neq 0.
\]

To obtain $\AAA9)$ observe that sets in $\col{I}$ are pairwise disjoint. Thus by Lemma~\ref{Hfunc} H7) we get
\[
    \sum_{j} \eLpp{\mardifold{g_I}{j}}{\Omega}{1}=\sum_{A\in\col{I}}\sum_{j} \eLpp{\mardifold{\khan{A}}{j}}{\Omega}{1} \leq 4\sum_{A\in\col{I}}\Pp(A)= 4 \Pp(\Est{\col{I}}).
\]

Now we turn to proving $\AAA10)$, starting with the right hand side of inequality \eqref{likedyad}. Clearly we have
\[
\Pp(\Est{\col{[0,1]}})= \Pp(\Io).
\]
We inductively assume that the fixed interval $I$ satisfies 
\begin{equation}\label{230901}
\Pp(\Io) \left(\frac{1}{2}-\varepsilon-\tilde{\varepsilon}\right)^{-\log_2(|I|)}\leq\Pp(\Est{\col{I}})\leq|I| \Pp(\Io).
\end{equation}
We prove that inequality \eqref{230901} holds for $I^{+}$ and $I^{-}$. We consider the estimates for $I^{+}$. By $\AAA3)$, $\AAA5)$ and Lemma \ref{Hfunc} H6) we have
\[
\Pp(\Est{\col{I^{+}}}) \leq \sum_{A\in\col{I}} \Pp(\{\khan{A}=1\})\leq \frac{1}{2}\sum_{A\in\col{I}} \Pp(A)=\frac{1}{2} \Pp(\Est{\col{I}})\leq\frac{1}{2}|I| \Pp(\Io).
\]
Since $|I^{+}|=\frac{1}{2}|I|$ we have verified the right hand side of inequality \eqref{230901}. We put $l=-\log_2|I|$. Observe that by Lemma~\ref{Hfunc} H6) and Lemma \ref{enhcond1} $iii)$ we have
\[
\Pp(\{\khan{A}=1\}\cap \Est{\Tc_{l+1}})\geq \left(\frac{1}{2}-\varepsilon-\tilde{\varepsilon}\right)\Pp(A)
\]
for any $A\in\Tc_{l}$. 
Therefore  
\[
\Pp(\Est{\col{I^{+}}})\stackrel{\AAA3)}{\geq} \sum_{A\in\col{I}}\Pp(\{\khan{A}=\xi\}\cap \Est{\Tc_{l+1}})\geq \left(\frac{1}{2}-\varepsilon-\tilde{\varepsilon}\right)\sum_{A\in\col{I}} \Pp(A)\stackrel{\eqref{230901}}{\geq} \left(\frac{1}{2}-\varepsilon-\tilde{\varepsilon}\right)^{l+1} \Pp(\Io)
\]
For any small enough $\hat{\varepsilon}>0$ we choose $\varepsilon>0$, $\tilde{\varepsilon}>0$ such that for every $q\in\{0,1,\ldots,n\}$ we have
\begin{equation}\label{230902}
\left(\frac{1}{2}-\varepsilon-\tilde{\varepsilon}\right)^{q}\geq \frac{1}{2^{q}}-\hat{\varepsilon}.    
\end{equation}
In view of \eqref{230902} this gives us the left hand side of \eqref{likedyad} in $\AAA10)$.
\end{proof}
\begin{proof}[\bf{Proof of Theorem} \ref{infincarl}]
Since $\crl{\Ec}=\infty$. By Lemma \ref{unicolAAA} the assumptions of Lemma \ref{AAAimpMar} are satisfied with constant $K=4$. Lemma \ref{AAAimpMar} implies that the Banach space $X$ is of Haar type p and by \eqref{estmartconst}
\[
HT_p(X)\leq 4^{\frac{2}{p}-1} T_p,
\]
where $T_p$ is defined in \eqref{typeX}.
\end{proof}

\section{Finite Carleson constant and upper $L^p$ estimates}\label{sec:fincarl}
Let $\FF_n$ be a purely atomic $\sigma$-algebra for $n\in\NN$. Recall that for the filtered probability space  $(\Omega,\, \FF,\, \left(\FF_n\right),\, \Pp)$ we have defined a collection $\Ec$ associated with it (see \eqref{defcole}). Recall also that the Carleson constant of $ \Ec$ is given by
\[
\crl{\Ec}= \sup_{I\in \Ec}\frac{1}{\Pp(I)}\sum_{\substack{J\subseteq I\\ J\in \Ec}}\Pp(J).
\]
In the last section we investigated the case where the Carleson constant of $\Ec$ is infinite. The main result there is Theorem \ref{infincarl}  which asserts that in case $\crl{\Ec}=\infty$, the inequality \eqref{typeX} implies non-trivial martingale type of the Banach space $X$. 

In the present section we consider filtrations $(\FF_n)$ for which the Carleson constant of $\Ec$ is finite. We will exploit the observation going back to Carleson and Garnett \cite{Carleson1975} that $\Ec$ can be decomposed as 
\[
\Ec= \Ec_{1}\cup\Ec_{2}\cup\Ec_{3}\cup\cdots\cup\Ec_{M},  \qquad M=\lfloor 4\crl{\Ec}+1\rfloor, 
\]
and for any $A\in \Ec_{j}$ the following inequality holds true
\begin{equation*}
\Pp(\Est{\GC{k}{\Ec_j}{A}})=\sum_{K\in \GC{k}{\Ec_j}{A}} \Pp(K)\leq 2^{-k} \Pp(A),\qquad\forall\; k\in\NN.
\end{equation*}

The main result of this section is Theorem \ref{thmfincarl} below which asserts that the inequality \eqref{typeX} is satisfied for any Banach space $X$ provided that $\crl{\Ec}<\infty$.
\begin{T1}\label{thmfincarl}
Let $(\Omega,\, \FF,\, \left(\FF_n\right),\, \Pp)$ be a filtered probability space, where each $\FF_n$ is purely atomic. Let $1<p\leq 2$ and $\crl{\Ec}<\infty$. For any Banach space $X$ and any $f\in L^p(\Omega,X)$ the following inequality holds:
\begin{equation}\label{alws}
\eLpx{f}\leq T_p\left(\chng{\eLpp{\EE(f|\FF_0)}{\Omega,\, X}{p}^p}+ \sum_{n=1}^{\infty}\eLpx{\mardif{f}{n}}^p   \right)^{\frac{1}{p}},     
\end{equation}
where
\[
T_p=\left(1+\frac{1}{1-2^{-\frac{1}{p}}}\right)2^{2-\frac{1}{p}}(1+MT_p(\RR)^{p})^{\frac{1}{p}} (4\crl{\Ec}+1)^{1-\frac{1}{p}}.
\]
\end{T1}
\begin{Rem}
The constant $T_p$ in the Theorem \ref{thmfincarl} is not optimal. It reflects our approach, which is based on the disjointification procedure of Carleson and Garnett.   
\end{Rem}

\subsection{Definitions and auxiliary results}
In this subsection we concentrate on special cases of inequality \eqref{alws}.We first determine an algebraic basis for the space of martingale differences associated to a single atom. Recall that  in opening paragraph of Section \ref{sec: not} we had modified the fields $\FF_n$ in such a way that every atom $A \in \AAA_n$ will be, either split up in strictly smaller atoms in $\AAA_{n+1}$, or else never split at any later step $m > n$.

  Fix $n\in\NN$ and $A\in\AAA_n$. Recall that $A=\bigcup_{j=1}^{\Nay{A}} A_j$, where $A_j\in\AAA_{n+1}$ are ordered in such a way that
\begin{equation*}
 \Pp(A_k)= \max_{j\geq k} \Pp(A_{j}).     
\end{equation*}
 For any $j\in\{2,\ldots, \Nay{A}\}$ we put
\begin{equation}\label{defbasis}
 \hax{A_j}{\omega}= \chs_{A_j}(\omega) - \frac{\Pp(A_j)}{\Pp(A_{j-1})} \chs_{A_{j-1}}(\omega),\qquad \omega\in\Omega.  
\end{equation}

The functions $\{k_{A_j}\}_{j=2}^{N(A)}$ form an algebraic basis for the space of martingale differences restricted to $A$, that is, for the finite dimensional space 
\[
\{\left(\EE(f\big|\FF_{n+1})-\EE(f\big|\FF_{n})\right)\big|_A\;:\;f\in L^1\}.
\]

We put $\ha{\Omega}=\chs_{\Omega}$.  Observe that by \eqref{defbasis} and the definition of $\Ec$ \eqref{defcole} for every atom $A\in \Ec$ the function $\ha{A}$ is well defined and that $\{k_{A_j}\;:\; j\in\{2, \ldots, \Nay{A}\},\; A\in\AAA\}\cup\{\chs_{\Omega}\}$ is just an alternative listing of $\{k_B\;:\;B\in\Ec \}.$ Note that $\Omega\in\Ec$ and that for $1\leq p < \infty$ we have
\[
\overline{\operatorname{span}\{k_B\;:\;B\in\Ec \}}^{L^p}= L^p(\Omega,\FF, \Pp).
\]
\begin{figure}[h]
\begin{center}
\includegraphics[width=5cm ]{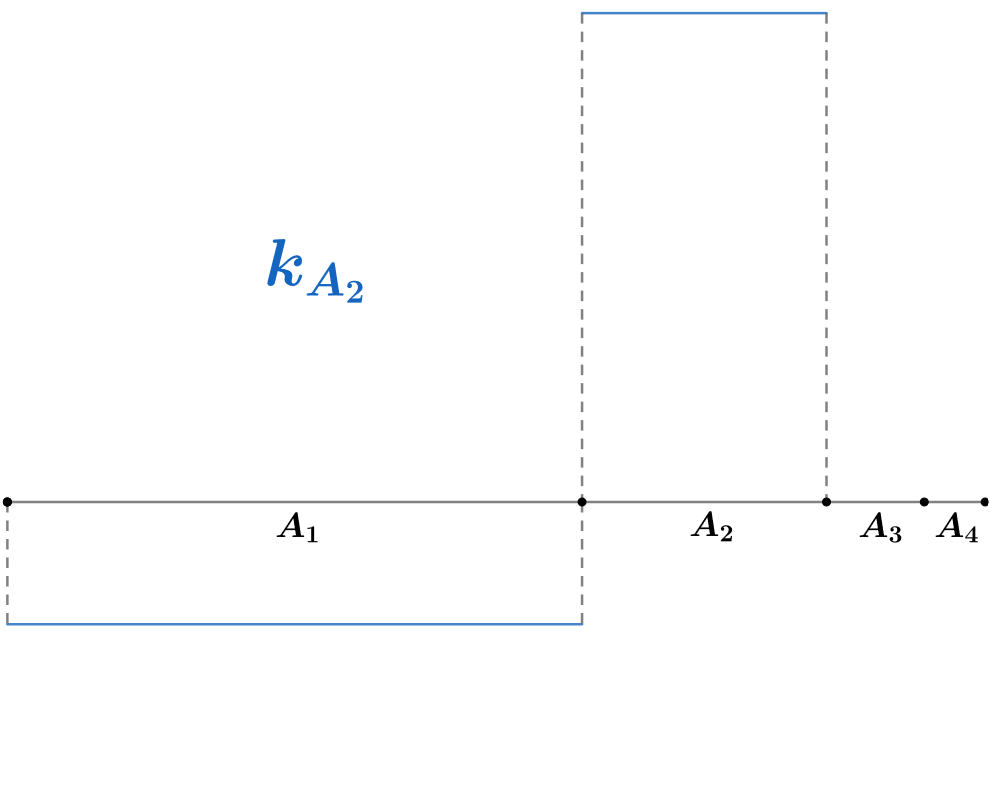}
\includegraphics[width=5cm ]{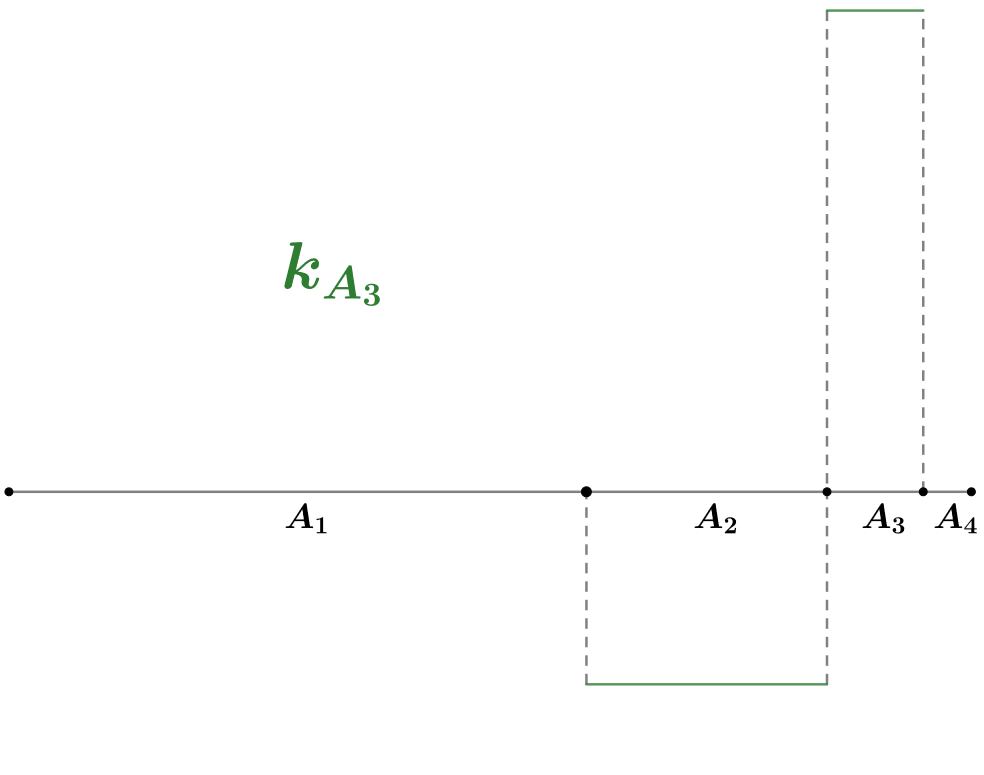}
\includegraphics[width=5cm ]{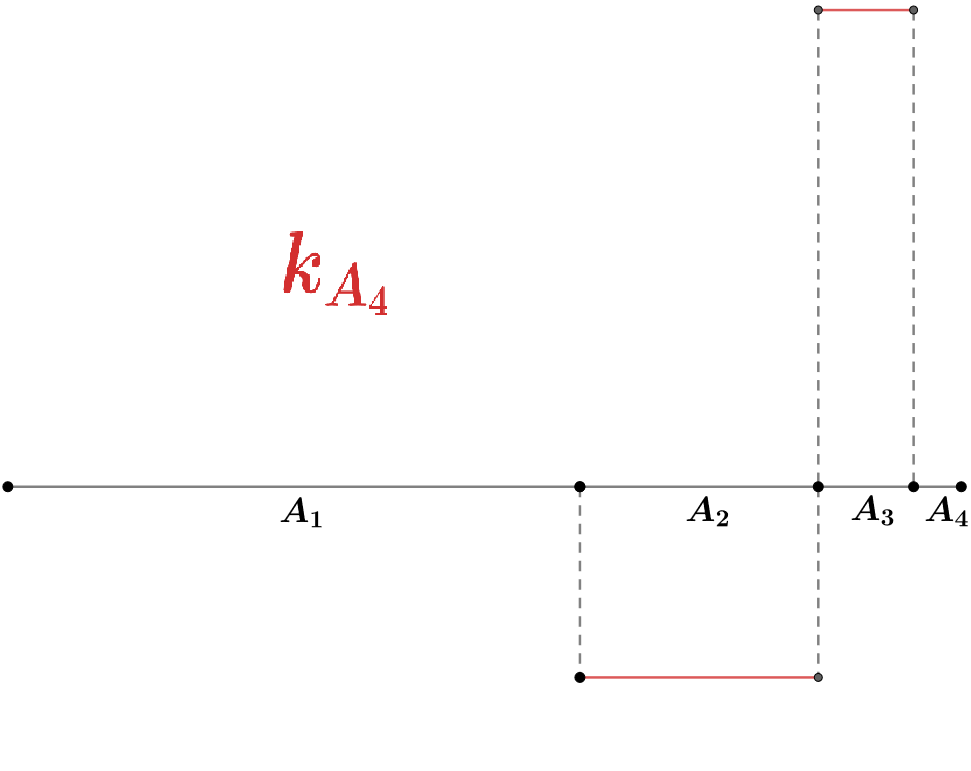}
\caption{Given $A\in\AAA_n$, the figure depicts $A_1,A_2,A_3,A_4\in\AAA_{n+1}$ such that $A= A_1\cup A_2\cup A_3\cup A_4$. We have $N(A)=4$, and the functions $\{\ha{A_2}$, $\ha{A_3}$, $\ha{A_4}\}$ (defined by \eqref{defbasis}) form an algebraic basis of the space of martingale differences restricted to $A$.}
\label{fig: wykreshax}
\end{center}
\end{figure}

Recall that for $A\in\AAA$ we defined $A_1$ by \eqref{notorder} and we put $\st{A}=A_1$. Now we let 
\begin{equation}\label{deftiAsun}
 \ti{A}:=\bigcup\limits_{j=2}^{N(A)}  A_j=A\backslash \st{A},
\end{equation}
put 
\begin{equation}\label{defrozhaa}
 \haxi{A}{\omega}=\Pp(\ti{A})\chs_{\st{A}}(\omega)  -  \Pp(\st{A})\chs_{\ti{A}}(\omega), \quad \omega\in\Omega ,
\end{equation}
and finally we define the collection
\begin{equation}\label{collectb}
 \bbb=\{\ti{A}: A\in \mathscr{A}\}.   
\end{equation}
\begin{figure}[h]
\begin{center}
\includegraphics[width=6.3cm ]{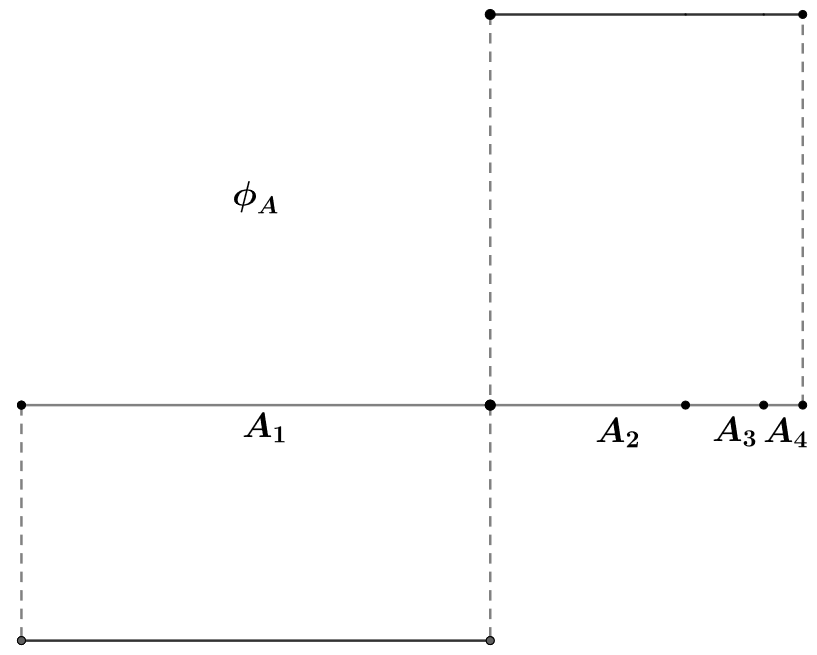}
\caption{Given $A\in\AAA_n$, the figure depicts $A_1,A_2,A_3,A_4\in\AAA_{n+1}$ such that $A= A_1\cup A_2\cup A_3\cup A_4$. Then $\ti{A}=A_2\cup A_3\cup A_4$ and $\st{A}=A_1$. The figure depicts the graph of $\hai{A}$.}
\label{fig: wykresphi}
\end{center}
\end{figure}
\begin{Rem}\label{rem :collectionBAC}
In the following list of remarks we comment on the relation between the collections $\AAA$, $\bbb$, $\ccc$ and $\Ec$.
\begin{enumerate}[label=\ref{rem :collectionBAC}.\arabic*), font=\normalfont]
    \item\label{rem :collectionBAC1} Given $K\in \bbb$ there exist (unique) $n\in\NN\cup\{0\}$ and $A(K)\in\AAA_n$ such that $K=\ti{A(K)}$. In that case there exists uniquely defined collection consisting of pairwise disjoint atoms $\{A_j\;:\; 2\leq j\leq N(A(K))\}\subset\AAA_{n+1}$ such that $K=\bigcup_{j=2}^{N(A(K))} A_j.$ Moreover $A_j\in\Ec$ for $j\in\{2,...,N(A(K))\}$. 
    \item\label{rem :collectionBAC2} If $A\in\AAA$ and $K\in\bbb$ satisfy $K\subsetneq \ti{A}$, then there exists $j\in\{2,\ldots, N(A)\}$ such that $K\subset A_j$. 
    \item\label{rem :collectionBAC3} If $A,B\in\AAA$ and there exist $j\in\{2,...,N(A)\}$ such that $A_j\subsetneq B$ then $\ti{A}\subsetneq B$. Similarly if $A_j\subseteq \ti{B}$ then $\ti{A}\subseteq \ti{B}.$
    \item\label{rem :collectionBAC4} Given $L\in \ccc$ there exist (unique) $n\in \NN\cup\{0\}$ and $A(L)\in\AAA_n$ such that $L=\st{A(L)}.$
    \item\label{rem :collectionBAC5} If $K\in\bbb$, $L\in\ccc$ and $A\in\AAA$ such that $A=A(K)$ and $A=A(L)$, then necessarily  $K=\ti{A}$ and $L=\st{A}$. 

\end{enumerate}
In summary by \eqref{deftiAsun}, $A\in\AAA$  is split into $\ti{A}\in\bbb$ and $\st{A}\in\ccc$; conversely by Remark \ref{rem :collectionBAC5}, $\ti{A}\in\bbb$ uniquely determines $A\in\AAA$ and $\st{A}\in\ccc$ ; and again by Remark \ref{rem :collectionBAC5}, $\st{A}\in\ccc$ uniquely determines $A\in\AAA$ and $\ti{A}\in\bbb.$ 
 
 \begin{equation}
\begin{tikzpicture}
  \matrix (m) [matrix of math nodes,row sep=2.5em,column sep=2.5em,minimum width=1em]
  {
      & A & & & A& & &A&  \\
     \ti{A}& &\st{A}\quad; &\ti{A}&  & \st{A}\quad; &\ti{A}& &\st{A}&\\
  };
  \path[-stealth]
   (m-1-2) edge node [above] {} (m-2-1)
    (m-1-2) edge  node [above] {} (m-2-3)
      (m-2-4) edge node [above] {} (m-1-5)
      (m-2-4)  edge node [above] {} (m-2-6)
        (m-2-9) edge node [above] {} (m-1-8)
      (m-2-9)  edge node [above] {} (m-2-7)
          ;
\end{tikzpicture}
\end{equation}

\end{Rem}

\begin{L1}\label{finBcEc}
The Carleson constant of the collection $\Ec$ is finite if and only if the Carleson constant of the collection $\bbb$ is finite. Moreover
\[
\crl{\bbb} \leq \crl{\Ec} \leq 1+\crl{\bbb} 
\]
\end{L1}
\begin{proof}
First we assume that $\crl{\Ec}<\infty$. Observe that for any $\ti{A}\in\bbb$ we have 
\[
\begin{split}
\sum_{\substack{\ti{B}\subseteq \ti{A}\\ \ti{B}\in \bbb}}\Pp(\ti{B})&=\Pp(\ti{A})+\sum_{\substack{\ti{B}\subsetneq \ti{A}\\ \ti{B}\in \bbb}}\Pp(\ti{B})
\\&\stackrel{Rem.\ref{rem :collectionBAC2}}{=}\Pp(\ti{A})+\sum_{j=2}^{\Nay{A}}\sum_{\substack{\ti{B}\subsetneq A_j\\ \ti{B}\in \bbb}}\Pp(\ti{B})
\\&\stackrel{Rem.\ref{rem :collectionBAC1}}{=}\sum_{j=2}^{\Nay{A}}\Pp(A_j)+\sum_{j=2}^{\Nay{A}}\sum_{\substack{\ti{B}\subsetneq A_j\\ \ti{B}\in \bbb}}\Pp(\ti{B})
\\&\stackrel{Rem.\ref{rem :collectionBAC1}}{=}\sum_{j=2}^{\Nay{A}}\Pp(A_j)+\sum_{j=2}^{\Nay{A}}\sum_{\substack{\ti{B}\subsetneq A_j\\ \ti{B}\in \bbb}}\sum_{k=2}^{\Nay{B}}\Pp(B_k)
\\ &=\sum_{j=2}^{\Nay{A}}\sum_{\substack{B \subseteq A_j\\ B\in \Ec}}\Pp(B).
\end{split}
\]
Hence by definition of $\crl{\Ec}$ we obtain
\[
\sum_{\substack{\ti{B}\subseteq \ti{A}\\ \ti{B}\in \bbb}}\Pp(\ti{B})\leq \crl{\Ec}\sum_{j=2}^{\Nay{A}} \Pp(A_j)=\crl{\Ec}\Pp(\ti{A}).
\]
Therefore
\[
\crl{\bbb}\leq \crl{\Ec} <\infty.
\]
On the other hand for $A\in \Ec$ and $\crl{\bbb}<\infty$ we have
\begin{equation}\label{eq: bcolzero}
\begin{split}
\sum_{\substack{B\subseteq A\\ B\in \Ec}}\Pp(B) &=\Pp(A)+\sum_{C\in \GC{1}{\Ec}{A}}\left(\Pp(C)+\sum_{\substack{B\subsetneq C\\ B\in \Ec}}\Pp(B) \right).
\end{split}
\end{equation}
For a given  $C \in \GC{1}{\Ec}{A}$ such that $C\in\AAA_n$ there exists a unique $D\in\AAA_{n-1}$ such that $C\subseteq \ti{D},$ and hence for any $B\subsetneq C$ we have $B \subsetneq \ti{D}$. Moreover by Remark \ref{rem :collectionBAC1} there exist $\{D_j\;:\; 2\leq j\leq N(D)\}\subset \AAA_n$ such that $\ti{D}=\bigcup_{j=2}^{N(D)} D_j$. There exists unique $m\in\{2,\ldots,N(D)\}$ such that $C=D_m$. By Remark \ref{rem :collectionBAC3} $C\subsetneq A$ implies $\ti{D}\subsetneq A$, and hence $D_j\subsetneq A$ for any $j\in\{2,\ldots,N(D)\}$. Moreover $\ti{D}\in \GC{1}{\bbb}{A}$ and $D_j\in \GC{1}{\Ec}{A}$ for any $j\in\{2,\ldots,N(D)\}$.   It follows
\begin{equation}\label{eq: bcolfirst}
\begin{split}
 \sum_{C\in \GC{1}{\Ec}{A}}\left(\Pp(C)+\sum_{\substack{B\subsetneq C\\ B\in \Ec}}\Pp(B) \right) &=\sum_{\ti{D}\in \GC{1}{\bbb}{A}}\left(\sum_{j=2}^{N(D)}\Pp(D_j)+\sum_{\substack{B\subsetneq D_j\\ B\in \Ec}}\Pp(B) \right)
\\&=\sum_{\ti{D}\in \GC{1}{\bbb}{A}}\left(\Pp(\ti{D})+\sum_{\substack{B\subsetneq \ti{D}\\ B\in \Ec}}\Pp(B) \right). \end{split}  
\end{equation}
Fix $\ti{D}\in\bbb$. Similarly by Remark \ref{rem :collectionBAC1} and Remark \ref{rem :collectionBAC3} we obtain
\begin{equation}\label{eq: bcolsecond}
P(\ti{D})+\sum_{\substack{B\subsetneq \ti{D}\\ B\in \Ec}}\Pp(B)=P(\ti{D})+\sum_{\substack{\ti{E}\subsetneq \ti{D}\\ \ti{E}\in \bbb}}\Pp(\ti{E})\stackrel{\eqref{defcrl}}{\leq}\crl{\bbb}\Pp(\ti{D}).
\end{equation}
Combining \eqref{eq: bcolzero}, \eqref{eq: bcolfirst} and \eqref{eq: bcolsecond} we get 
\[
\begin{split}
\sum_{\substack{B\subseteq A\\ B\in \Ec}}\Pp(B)&\leq \Pp(A)+\crl{\bbb}\sum_{\ti{D}\in \GC{1}{\bbb}{A}}\Pp(\ti{D})
\\ &\leq (1+\crl{\bbb}) \Pp(A).
\end{split}
\]
Hence
\[
\crl{\Ec} \leq 1+\crl{\bbb} <\infty.
\]
\end{proof}
Now we consider nested collections with finite Carleson constants and we review the disjointification lemma originally due to Carleson and Garnett \cite{Carleson1975}. (Also see \cite{Mueller2005}.)
\begin{L1}\label{smmes}
Let $\Rc$ be a nested collection of measurable sets such that $\crl{\Rc}<\infty$ and $M_{\Rc}=\lfloor 4\crl{\Rc}+1 \rfloor$. There are families $\{\Rc_j\}^{M_{\Rc}-1}_{j=0}\subseteq \Rc$ such that
\[
\Rc = \bigcup_{l=0}^{M_{\Rc}-1} \Rc_l,\quad\mbox{and}\quad  \Rc_j\cap \Rc_i =\emptyset,
\]
where $0\leq i< j \chng{ < M_{\Rc}}$. Moreover
\[
\sum_{J\in \GC{1}{\Rc_{j}}{A}} \Pp(J)\leq \frac{\Pp(A)}{2}, \qquad\forall\; A\in \Rc_j\qquad\forall\;j\in\{0,\cdots, M_{\Rc}-1\}.
\]
\end{L1}
In consequence we obtain a bound on the size of subsequent generations
\begin{equation}\label{geomet}
\sum_{K\in \GC{k}{\Rc_i}{A}} \Pp(K)=\sum_{J\in \GC{k-1}{\Rc_i}{A}} \sum_{K\in \GC{1}{\Rc_i}{J}} \Pp(K)\leq \frac{1}{2} \sum_{J\in \GC{k-1}{\Rc_i}{A}} \Pp(J)\leq 2^{-k} \Pp(A).
\end{equation}
Equation \eqref{geomet} asserts that the sequence $\Pp(\Est{\GC{k}{\Ec_i}{A}})\Pp(A)^{-1}$ decreases at geometric rate uniformly with respect to $A\in\Rc$. As a result we will be able to prove that $$\{\chs_{K}\}_{K\in\Ec_i}$$ forms an almost disjointly supported sequence of functions. See Lemma \ref{conscol} and Lemma \ref{estcolb} below. Here we say that the sequence of functions is almost disjointly supported if it satisfies conditions \eqref{nosnikgj} and \eqref{malenie} below.  Lemma \ref{estDD} links almost disjointly supported sequences of functions to martingale type estimates.
\begin{L1}\label{estDD}
Let $X$ be a Banach space, $g_j : (\Omega,\Sigma, \Pp) \rightarrow X$, $j\in\NN$ denote a sequence of measurable functions and $1\leq p <\infty$. Assume that there exists a sequence $\{a_k\}_{k\in \NN}\in \ell^1(\NN)$ and a family of measurable sets $\{\DD{k}\}_{k\in\NN}$  such that 
\begin{equation}\label{nosnikgj}
   \Pp\left(\supp{g_j} \;\backslash\;  \bigcup_{k\geq j} \DD{k}\right)=0, 
\end{equation}
and, for $k,j \in\NN$,
\begin{equation}\label{malenie}
  \int_{\DD{k+j-1}} \nx{g_{j}}^p\dP \leq |a_{k}|^p \eLpx{g_j}^p,  
\end{equation}
then
\[
\eLpxb{\sum_{j\in\NN} g_j} \leq \bigg(  \sum_{k\in\NN} |a_k|\bigg) \left(\sum_{j\in\NN} \eLpx{g_j}^p\right)^{\frac{1}{p}}.
\]
\end{L1}
\begin{proof}
The claim follows from Minkowski's inequality. Let $E=\bigcap_{k=1}^{\infty}\DD{k};$ it follows from \eqref{malenie} that
\[
\int_E \nx{g_j}^pd\Pp=0\quad\mbox{for every}\quad j\geq 1,
\]
hence if $F=\bigcup_{k=1}^{\infty}\DD{k}$ and $U_k= \DD{k}\backslash \bigcup_{l>k} \DD{l}$ then $F\backslash E=\bigcup_{k=1}^{\infty} U_k$ and
\[
\begin{split}
\bigg(\int_{\Omega}\nxb{\sum_{j\in\NN} g_j}^p\dP\bigg)^{\frac{1}{p}}&=\bigg(\int_{F\backslash E}\nxb{\sum_{j\in\NN} g_j}^p\dP\bigg)^{\frac{1}{p}}\\&=\bigg(\sum_{k\in\NN}\int_{U_k}\nxb{\sum_{j\in\NN} g_j}^p\dP\bigg)^{\frac{1}{p}} 
\stackrel{\eqref{nosnikgj}}{=}\bigg(\sum_{k\in\NN}\int_{U_k}\nxb{\sum_{j\in\NN\cup\{0\}} \chs_{j< k} g_{k-j}}^p\dP\bigg)^{\frac{1}{p}}.
\end{split}
\]
Therefore
\[
\begin{split}
\bigg(\int_{\Omega}\nxb{\sum_{j\in\NN} g_j}^p\dP\bigg)^{\frac{1}{p}}  &\leq\bigg(\sum_{k\in\NN}\int_{\DD{k}}\nxb{\sum_{j\in\NN\cup\{0\}} \chng{\chs_{j< k}} g_{k-j}}^p\dP\bigg)^{\frac{1}{p}} 
\\&\chng{=\bigg(\sum_{k\in\NN}\int_{\Omega} \chs_{\DD{k}}\nxb{\sum_{j\in\NN\cup\{0\}} \chs_{j< k} g_{k-j}}^p\dP\bigg)^{\frac{1}{p}}}
\\&\chng{=\bigg(\sum_{k\in\NN}\int_{\Omega} \nxb{\sum_{j\in\NN\cup\{0\}} \chs_{\DD{k}} \chs_{j< k} g_{k-j}}^p\dP\bigg)^{\frac{1}{p}}}
\\&\chng{=\bigg( \int_{\Omega} \sum_{k\in\NN} \nxb{\sum_{j\in\NN\cup\{0\}} \chs_{\DD{k}} \chs_{j< k} g_{k-j}}^p\dP\bigg)^{\frac{1}{p}}}
\end{split}
\]
 \chng{The last term we interpret as the $L^p(\ell^p(X))$ norm of a sum of vectors and analogously we interpret the following term 
\[
\sum_{j\in\NN\cup\{0\}}\bigg(\sum_{k\in\NN}\chng{\chs_{j< k}}\int_{\DD{k}}\nx{ g_{k-j}}^p\dP \bigg)^{\frac{1}{p}}
\]
 as the corresponding sum of norms in the $L^p(\ell^p(X))$. Therefore Minkowski's inequality in $L^p(\ell^p(X))$ gives
 \[
\bigg( \int_{\Omega} \sum_{k\in\NN} \nxb{\sum_{j\in\NN\cup\{0\}} \chs_{\DD{k}} \chs_{j< k} g_{k-j}}^p\dP\bigg)^{\frac{1}{p}} \leq\sum_{j\in\NN\cup\{0\}}\bigg(\sum_{k\in\NN}\chs_{j< k}\int_{\DD{k}}\nx{ g_{k-j}}^p\dP \bigg)^{\frac{1}{p}}.
 \]
 Summing up we obtain
\[
\begin{split}
\bigg(\int_{\Omega}\nxb{\sum_{j\in\NN} g_j}^p\dP\bigg)^{\frac{1}{p}} & \leq\sum_{j\in\NN\cup\{0\}}\bigg(\sum_{k\in\NN}\chng{\chs_{j< k}}\int_{\DD{k}}\nx{ g_{k-j}}^p\dP \bigg)^{\frac{1}{p}} 
\\&\stackrel{\eqref{malenie}}{\leq}\sum_{j\in\NN\cup\{0\}}\bigg(\sum_{k=j+1}^{\infty} |a_{j+1}|^p \eLpx{g_{k-j}}^p  \bigg)^{\frac{1}{p}} 
\\&= \bigg(\sum_{j\in\NN} |a_j|\bigg) \bigg(\sum_{n\in\NN} \eLpx{g_n}^p\bigg)^{\frac{1}{p}}.
\end{split}
\] }

\end{proof}
Using the above lemma we obtain an estimate on the norm of the function supported on the atoms from $\Ec$. 
\begin{L1}\label{conscol}
Let $\Rc$ be a nested collection of measurable subsets of $\Omega$ such that $\crl{\Rc} < \infty$. For any Banach space $X$ and any $\{x_A\}_{A\in\Rc} \subseteq X$ the following inequality holds true
\[
\eLpx{\sum\limits_{A\in\Rc} x_A \chs_{A}}\leq\stlemcrl{\Rc} \left(\sum_{A\in \Rc} \nx{x_A}^p \Pp(A)\right)^{\frac{1}{p}}  .
\]
\end{L1}
\begin{proof}
Let 
\[
f(\omega)=\sum\limits_{A\in\Rc} x_A \chs_{A}(\omega).
\]
We rewrite $f$ using families $\Rc_i$ from Lemma \ref{smmes}
\[
f(\omega)=\sum_{i=0}^{M_{\Rc}-1}\sum\limits_{A\in\Rc_i} x_A \chs_{A}(\omega).
\]
Recall that $M_{\Rc}=\lfloor 4 \crl{\Rc} +1 \rfloor$. We put
\[
f_i(\omega)=\sum\limits_{A\in\Rc_i} x_A \chs_{A}(\omega).
\]
It is enough to prove
\begin{equation}\label{forjin}
\eLpx{f_i}^p \leq C_p^p \sum\limits_{A\in\Rc_i}\nx{x_A}^p \Pp(A).
\end{equation}
Indeed by triangle and H\"older inequalities we would have
\[
\eLpx{f}^p\leq M_{\Rc}^{p-1} \sum_{j=0}^{M_{\Rc}-1}\eLpx{f_j}^p
\leq  C^p_p M_{\Rc}^{p-1}\sum_{j=0}^{M_{\Rc}-1} \sum\limits_{A\in\Rc_j} \nx{x_A}^p \Pp(A).
\]
Let $0\leq i,j\leq M_{\Rc}-1$. If $i\neq j$ then $\Rc_i\cap\Rc_j=\emptyset$. Moreover $\bigcup_{j=0}^{M_{\Rc}-1} \Rc_j = \Rc$. Hence
\[
\eLpx{f}\leq C_p M_{\Rc}^{1-\frac{1}{p}}\left( \sum_{A\in \Rc} \nx{x_A}^p \Pp(A)\right)^{\frac{1}{p}}.
\]
We return to the proof of \eqref{forjin}. We fix $i\in\{0,\ldots,\, M_{\Rc}-1\}$. First we assume that $\Omega\not\in\Rc_i$. For every $j\in\NN$ we put 
\[
g_j=\sum_{A \in \GC{j}{\Rc_i}{\Omega}} x_A \chs_A \qquad\mbox{and}\qquad \DD{j}= \Est{\GC{j}{\Rc_i}{\Omega}}.
\]
The sets $D_j=\Est{\GC{j}{\Rc_i}{\Omega}}$ are decreasing with respect to inclusion i.e. 
\[
\Omega\supset\DD{1}\supset\cdots\supset \DD{j}\supset\DD{j+1}\supset \cdots.
\] 
Therefore for $j\geq 1$ we have
\begin{equation}\label{eq: suppgjcollr}
\supp{g_j}=\DD{j}=\bigcup_{k\geq j} \DD{k}.
\end{equation}
Note that for any $A\in\GC{j}{\Rc_i}{\Omega}$ the function $g_j\big|_{A}$ is equal to the constant vector value $x_A$. Since $\Est{\GC{k}{\Rc_i}{A}}\subseteq A$ we have
\[
\frac{1}{\Pp(\Est{\GC{k}{\Rc_i}{A}})}\int_{\Est{\GC{k}{\Rc_i}{A}}} \nx{g_j}^p\dP=\frac{1}{\Pp(A)}\int_{A}\nx{g_j}^p\dP=\nx{x_A}^p
\]
Therefore
\begin{equation}\label{eq:geomtonatom}
\begin{split}
\int_{\Est{\GC{k}{\Rc_i}{A}}} \nx{g_j}^p\dP&=\frac{\Pp(\Est{\GC{k}{\Rc_i}{A}})}{\Pp(A)}\int_{A}\nx{g_j}^p\dP
\\&\stackrel{\eqref{geomet}}{\leq} 2^{-k}\int_{A}\nx{g_j}^p\dP.
\end{split}    
\end{equation}
Since $\{\GC{k}{\Rc_i}{A}\}_{A\in\GC{j}{\Rc_i}{\Omega}}$ is a family of pairwise disjoint sets and 
\[
\Est{\GC{k+j}{\Rc_i}{\Omega}}=\bigcup_{A\in\GC{j}{\Rc_i}{\Omega}} \Est{\GC{k}{\Rc_i}{A}}
\]
the estimates \eqref{eq:geomtonatom} yield
\begin{equation}\label{smallgeoint}
\begin{split}
\int_{\Est{\GC{k+j}{\Rc_i}{\Omega}}}\nx{g_j}^p\dP&=\sum_{A\in\GC{j}{\Rc_i}{\Omega}}\int_{\Est{\GC{k}{\Rc_i}{A}}} \nx{g_j}^p\dP
\\&\leq 2^{-k}\sum_{A\in\GC{j}{\Rc_i}{\Omega}}\int_{A}\nx{g_j}^p\dP
\\&=2^{-k}\int_{\Omega}\nx{g_j}^p\dP
\end{split}    
\end{equation}
Recall that $D_{k+j}=\Est{\GC{k+j}{\Rc_i}{\Omega}}$. Therefore \eqref{smallgeoint} and \eqref{eq: suppgjcollr} yield assumptions of Lemma \ref{estDD} with $g_j$, $\DD{j}$ and $a_k=2^{-\frac{k-1}{p}}$ for $k\in\NN$. Thus
\[
\eLpx{f_i}\leq C_p \left(\sum_{j=1}^{\infty}\eLpx{g_j}^p\right)^{\frac{1}{p}},
\]
where $C_p=(1-2^{-\frac{1}{p}})^{-1}$.
Since for every $j\in\NN$ the sets in $\GC{j}{\Rc_i}{\Omega}$ are pairwise disjoint we get 
\[
\eLpx{f_i}^p\leq C^p_p \sum_{j=1}^{\infty}\sum_{A \in \GC{j}{\Rc_i}{\Omega}} \nx{x_A}^p \Pp(A)=C^p_p\sum_{A \in \Rc_i}\nx{x_A}^p \Pp(A).
\]
Recall that in the above argument we have assumed that $\Omega\notin \Rc_i$. If $\Omega\in\Rc_i$ it suffices to change the definitions of $g_j$ and $\DD{j}$ as follows:
\[
g_0=x_{\Omega}\chs_{\Omega},\qquad \DD{0}= \Omega,
\]
\[
g_{j+1}=\sum_{A \in \GC{j}{\Rc_i}{\Omega}} x_A \chs_A,\qquad \DD{j+1}= \Est{\GC{j}{\Rc_i}{\Omega}}\qquad\mbox{where } j\in\NN,
\]
and repeat the above argument using these changes.
\end{proof}
The above lemma gives an estimate for the part of the function which is "contained" in the collection $\Ec$. We have similar bounds for the part "outside of the collection" i.e. linear span of the functions $\{\chs_{\st{A}}\}_{A\in\AAA}$.

Recall that for the collection $\Ec$ we have defined the corresponding collection $\bbb$ in \eqref{collectb}.

\begin{L1}\label{singlecol}
Let $1 < p \leq2$. Let $\Rc\subset\bbb$ and X be a Banach space. For any choice of $n\in\NN\cup\{0\}$,
\[
A\in \FF_{n}\quad \mbox{and}\quad \{x_{\ti{B}}\}_{\ti{B}\in \GC{1}{\Rc}{A}}\subset X
\]
we have 
\begin{equation*}
   \eLpx{\sum_{\ti{B}\in \GC{1}{\Rc}{A}} x_{\ti{B}}\Pp(\ti{B})\chs_{\st{B}} }\leq \stlem \left(\sum_{\ti{B}\in \GC{1}{\Rc}{A}} \nx{x_{\ti{B}}}^p \eLp{\hai{B}}^p\right)^{\frac{1}{p}},
\end{equation*}
where 
\[
\hai{B}=\Pp(\ti{B})\chs_{\st{B}}  -  \Pp(\st{B})\chs_{\ti{B}}.
\]

\end{L1}
\begin{proof}
We put
\[
f(\omega)=\sum_{\ti{B}\in \GC{1}{\Rc}{A}} x_{\ti{B}}\Pp(\ti{B})\chs_{\st{B}}(\omega)
\]
For any $A\in\Rc$ we define $T_A = A \backslash  \Est{\GC{1}{\Rc}{A}}$. See Figure \ref{fig: TI}.
\begin{figure}[t]
\begin{center}
\includegraphics[width=0.7\textwidth]{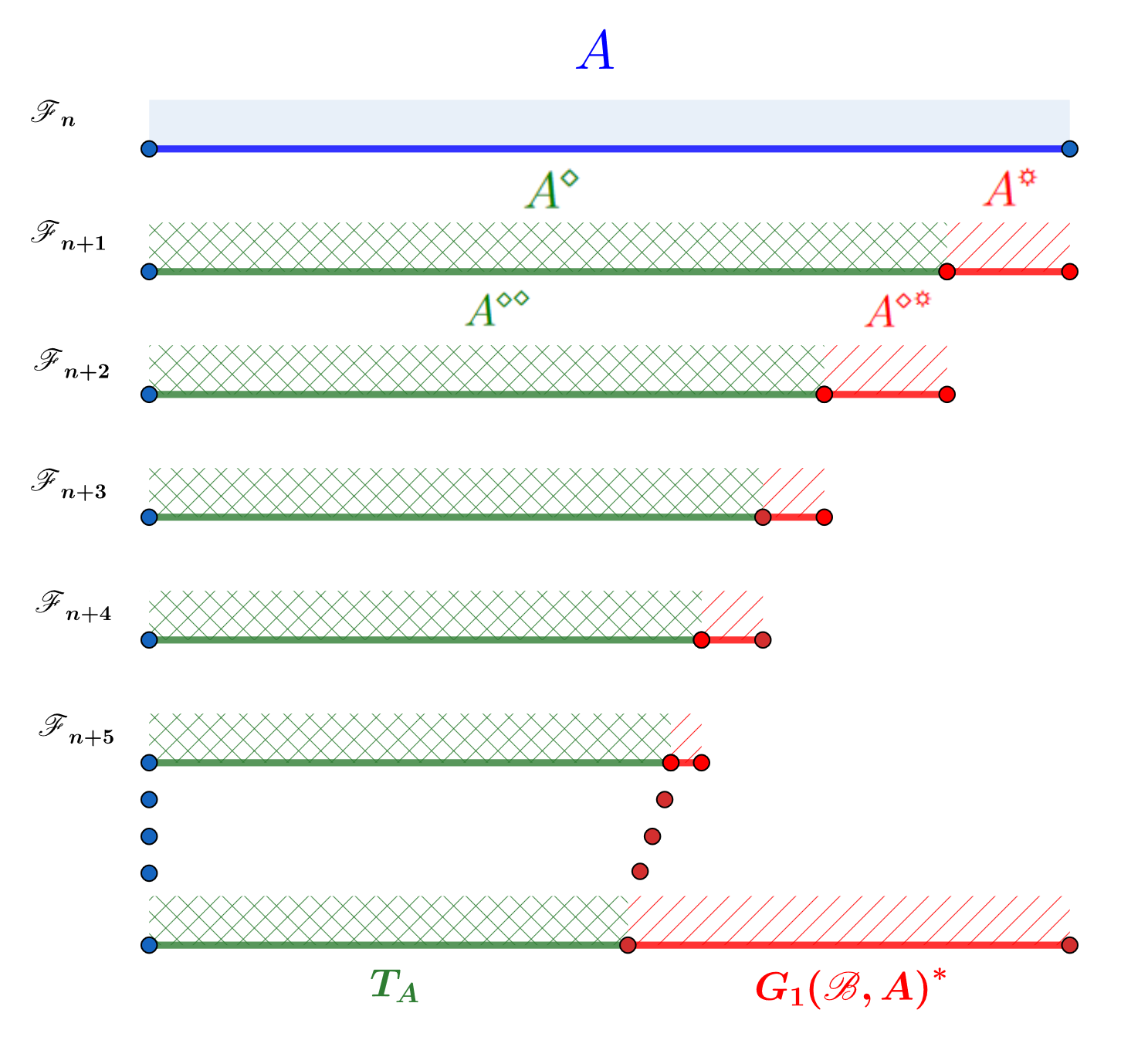}   
\caption{The picture highlights the different role of the atoms in $\AAA$. Let $A\in\AAA_n$ and $\Rc=\bbb$.  The pairwise disjoint red coloured intervals depict $\ti{A}$,$\ti{\st{A}}$, $\ti{\st{\st{A}}}$ etc. They form $\GC{1}{\Rc}{A}$. The intersection of the decreasing green coloured intervals forms the set $T_A$. }
\label{fig: TI}
\end{center}

\end{figure}
Since $\GC{1}{\Rc}{A}$ is a family of disjoint sets we have 
\[
\int_{\Omega} \nx{f}^p \dP=\int_{A} \nx{f}^p \dP= \int_{T_A} \nx{f}^p \dP + \sum_{\ti{B}\in \GC{1}{\Rc}{A}}\int_{\ti{B}} \nx{f}^p  \dP.
\]
The function $f$ is constant on every set $\ti{B}\in\GC{1}{\Rc}{A}$. For any $\ti{B}\in\GC{1}{\Rc}{A}$ by definition \eqref{defrozhaa} we have 
\begin{equation}\label{eq: TAIND}
\hai{B}\big|_{T_A}\equiv \Pp(\ti{B})\chs_{\st{B}}\big|_{T_A}
\end{equation}
Hence we have 
\begin{equation}\label{eq: TApart}
\begin{split}
\int_{T_A}\nx{f}^p \dP&= \int_{T_A}\nx{ \sum_{\ti{B}\in\GC{1}{\Rc}{A}} x_{\ti{B}}\Pp(\ti{B})\chs_{\st{B}}(\omega)}^p \dPx{\omega}
\\ &\leq \int_{T_A}\left| \sum_{\ti{B}\in\GC{1}{\Rc}{A}} \Pp(\ti{B})\chs_{\st{B}}(\omega) \nx{x_{\ti{B}}} \right|^p \dPx{\omega}
\\&\stackrel{\eqref{eq: TAIND}}{=}\int_{T_A}\left|\sum_{\ti{B}\in\GC{1}{\Rc}{A}}\haxi{B}{\omega} \nx{x_{\ti{B}}}\right|^p \dPx{\omega}.
\end{split}
\end{equation}
Now we fix a set $\ti{B}\in\GC{1}{\Rc}{A}$. 
Note that for $\omega\in\ti{B}$ we have
\[
\begin{split}
\nx{f(\omega)}^p 
&=\nxb{\sum_{\ti{K}\in\GC{1}{\Rc}{A}} x_{\ti{K}}\Pp(\ti{K})\chs_{\st{K}}(\omega)}^p\leq \bigg( \sum_{\ti{K}\in\GC{1}{\Rc}{A}} \nx{x_{\ti{K}}}\Pp(\ti{K})\chs_{\st{K}}(\omega)\bigg)^p
\\&=\bigg(-\nx{x_{\ti{B}}} \haxi{B}{\omega}+\nx{x_{\ti{B}}}\haxi{B}{\omega} + \sum_{\ti{K}\in\GC{1}{\Rc}{A}}\nx{x_{\ti{K}}}\Pp(\ti{K})\chs_{\st{K}}(\omega)\bigg)^p
\end{split}
\]
Since $\ti{B}\in\GC{1}{\Rc}{A}$ we get $\hai{K}(\omega)=\Pp(\ti{K})\chs_{\st{K}}(\omega)$ for $\ti{K}\in\GC{1}{\Rc}{A}$, $\ti{K}\neq\ti{B}$ and $\omega\in\ti{B}$  (note that $\chs_{\ti{K}}(\omega)=0$ on $\ti{B}$). Thus for $\omega\in\ti{B}$ we have
\[
\begin{split}
\nx{f(\omega)}^p 
&\leq \bigg| -\nx{x_{\ti{B}}}\haxi{B}{\omega}+\nx{x_{\ti{B}}}\haxi{B}{\omega} + \sum_{\ti{K}\in\GC{1}{\Rc}{A}\backslash \ti{B}}\nx{x_{\ti{K}}}\haxi{K}{\omega}\bigg|^p
\\&\leq 2^{p-1} \nx{x_{\ti{B}}}^p|\haxi{B}{\omega}|^p+2^{p-1} \bigg|\sum_{\ti{K}\in\GC{1}{\Rc}{A}}\haxi{K}{\omega}\nx{x_{\ti{K}}} \bigg|^p.
\end{split}
\]
Integrating over $\ti{B}$ gives us
\begin{equation}\label{eq: bsunpart}
\begin{split}
\int_{\ti{B}}\nx{f}^p \dP
&\leq 2^{p-1}\int_{\ti{B}} \nx{x_{\ti{B}}}^p|\hai{B}|^p\dP+2^{p-1}\int_{\ti{B}}   \bigg|\sum_{\ti{K}\in\GC{1}{\Rc}{A}}\hai{K} \nx{x_{\ti{K}}} \bigg|^p\dP
\\&\leq 2^{p-1}\nx{x_{\ti{B}}}^p\eLp{ \hai{B}}^p+2^{p-1}\int_{\ti{B}}   \bigg|\sum_{\ti{K}\in\GC{1}{\Rc}{A}}\hai{K} \nx{x_{\ti{K}}} \bigg|^p\dP.
\end{split}
\end{equation}
We sum inequalities \eqref{eq: bsunpart} over $\ti{B}\in\GC{1}{\Rc}{A}$ and we add \eqref{eq: TApart} obtaining
\begin{equation}\label{calosc}
\int_{\Omega} \nx{f}^p \dP\leq 2^{p-1} \sum_{\ti{B}\in\GC{1}{\Rc}{A}} \nx{x_{\ti{B}}}^p\eLp{\hai{B}}^p+ 2^{p-1}\int_{\Omega}\bigg|\sum_{\ti{B}\in\GC{1}{\Rc}{A}}\hai{B} \nx{x_{\ti{B}}} \bigg|^p\dP.
\end{equation}
We define
\[
g(\omega)= \sum_{\ti{B}\in\GC{1}{\Rc}{A}}\nx{x_{\ti{B}}}\haxi{B}{\omega}.
\]
By \eqref{typescalar} the Banach space $\RR$ has Martingale type $p$ with constant $MT_p(\RR)$. Hence
\begin{equation}\label{calka}
\begin{split}
 \int_{\Omega}\bigg|\sum_{\ti{B}\in\GC{1}{\Rc}{A}}\hai{B} \nx{x_{\ti{B}}} \bigg|^p\dP 
 &\leq MT_p(\RR)^{p}\sum_{n=1}^{\infty} \eLpx{\mardif{g}{n}}^p
 \\& =MT_p(\RR)^{p} \sum_{\ti{B}\in\GC{1}{\Rc}{A}}\eLp{\hai{B}}^p\nx{x_{\ti{B}}}^p.   
\end{split}
\end{equation}
By \eqref{calosc} and \eqref{calka} we have
\[
\int_{\Omega} \nx{f}^p \dP\leq 2^{p-1}(1+MT_p(\RR)^{p})\sum_{\ti{B}\in\GC{1}{\Rc}{A}}\eLp{\hai{B}}^p\nx{x_{\ti{B}}}^p.
\]
\end{proof}

The next Lemma will combine Lemma \ref{estDD} and Lemma \ref{singlecol} to obtain upper bounds on the $L^p$ norm of functions $f\in\operatorname{span}\{\chs_{\st{A}}: A\in\AAA \}$. We associated to the collection $\AAA$ collections $\bbb$ and $\ccc$ defined respectively in \eqref{collectb} and \eqref{defcole}. Recall that by Remark \ref{rem :collectionBAC} we have 
\[\{\st{A}: A\in\AAA \}=\{\st{A}:\ti{A}\in\bbb\}.\]
By Lemma \ref{finBcEc}, $\crl{\bbb}\leq\crl{\Ec}<\infty$. Thus by Lemma \ref{smmes} we can decompose the collection $\bbb$ into collections $\bbb_i$ for $i\in\{0, 1,2, \ldots , M_{\bbb}-1 \}$ such that
\begin{equation}\label{eq: decreaseofB_i}
    \sum_{B\in \GC{1}{\bbb_i}{A}} \Pp(B)\leq \frac{\Pp(A)}{2}, \qquad A\in \bbb_i,
\end{equation}
where $M_{\bbb}\leq 4\crl{\Ec}+1$.

\begin{L1}\label{estcolb}
Let $1<p\leq 2$  and $\crl{\Ec}<\infty$. For any Banach space $X$ and $\{x_{\ti{A}}\}_{\ti{A}\in\bbb}\subset X$ the function
\[
f(\omega)= \sum_{\ti{A}\in \bbb} x_{\ti{A}}\Pp(\ti{A})\chs_{\st{A}}(\omega)
\]
satisfies
\begin{equation}\label{outcolest}
\eLpx{f}\leq \tilde{C}_p \left(\sum_{\ti{A}\in \bbb} \nx{x_{\ti{A}}}^{p}. \eLp{\hai{A}}^p\right)^{\frac{1}{p}},   
\end{equation}
where
\[
\tilde{C}_p=\left(1+\frac{1}{1-2^{-\frac{1}{p}}}\right)\stlem (4\crl{\Ec}+1)^{1-\frac{1}{p}}.
\]

\end{L1}
\begin{proof}
Decompose $\bbb$ as $\bbb=\bbb_0\cup \bbb_1\cup\ldots\cup\bbb_{M_{\bbb}-1}$ according to Lemma \ref{smmes}. We define 
\[
f_i=\sum_{\ti{A}\in \bbb_i} x_{\ti{A}}\Pp(\ti{A})\chs_{\st{A}}.
\]
As in Lemma \ref{conscol} we get 
\[
\eLpx{f}^p\leq M_{\bbb}^{p-1} \sum_{i=0}^{M_{\bbb}-1} \eLpx{f_i}^p. 
\]
Since $\bbb_i$ are pairwise disjoint families of sets it is enough to obtain the estimate \eqref{outcolest} for the function $f_i$ and $i\in\{0,1,\ldots ,M_{\bbb} -1\}$. Clearly
\[
f_i=\sum_{j=1}^{\infty} g_j,
\]
where
\[
g_j:=\sum_{\ti{A}\in \GC{j}{\bbb_i}{\Omega}}  x_{\ti{A}}\Pp(\ti{A})\chs_{\st{A}}.
\]
Observe that for $j>1$
\[
\begin{split}
\int_{\Omega} \nx{g_j}^p\dP
&=\sum_{\ti{B}\in \GC{j-1}{\bbb_i}{\Omega}}\int_{\ti{B}} \nxb{\sum_{\ti{A}\in \GC{1}{\bbb_i}{\ti{B}}}  x_{\ti{A}}\Pp(\ti{A})\chs_{\st{A}}}^p\dP
\\&=\sum_{\ti{B}\in \GC{j-1}{\bbb_i}{\Omega}}\int_{\Omega} \nxb{\sum_{\ti{A}\in \GC{1}{\bbb_i}{\ti{B}}}  x_{\ti{A}}\Pp(\ti{A})\chs_{\st{A}}}^p\dP
\end{split}
\]
By Lemma \ref{singlecol} applied to $\Rc=\bbb_i$ and $\ti{B}\in \GC{j-1}{\bbb_i}{\Omega}$ we get
\begin{equation}\label{highordcol}
 \begin{split}
\int_{\Omega} \nx{g_j}^p\dP&\leq \stlemp \sum_{\ti{B}\in \GC{j-1}{\bbb_i}{\Omega}}\;\sum_{\ti{A}\in \GC{1}{\bbb_i}{\ti{B}}} \nx{x_{\ti{A}}}^p\eLp{\hai{A}}^p 
\\&= \stlemp\sum_{\ti{A}\in \GC{j}{\bbb_i}{\Omega}} \nx{x_{\ti{A}}}^p\eLp{\hai{A}}^p.
\end{split}   
\end{equation}
Observe that $\supp{g_{1}}\subset\Omega$ and for $j\geq 2$ we know that $\supp{g_{j}}\subset \GC{j-1}{\bbb_i}{\Omega}^*$ and that the function $g_{j}$ is constant on every set $\ti{A}\in\GC{j+1}{\bbb_i}{\Omega}$ (see Figure \ref{fig:constsupp}).

\begin{figure}[h!]
    \includegraphics[width=\textwidth]{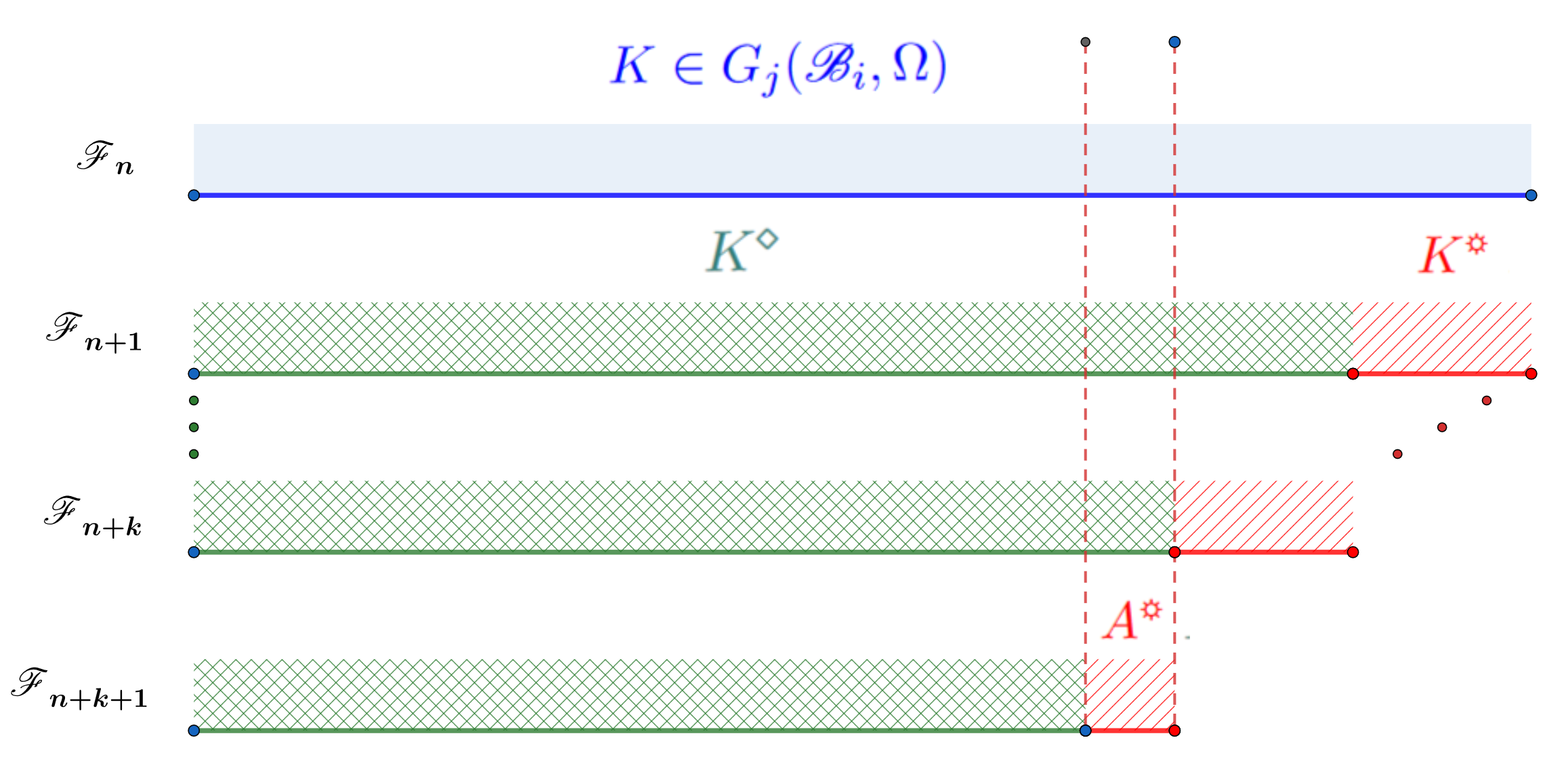}
    \caption{ This picture illustrates why functions $g_j$ have constant values on the sets from family $\GC{j+1}{\bbb_i}{\Omega}$. For $\ti{A}\in\GC{j+1}{\bbb_i}{\Omega}$ there exists a unique $K=\ti{B}\in\GC{j}{\bbb_i}{\Omega}$ such that $\ti{A}\in\GC{1}{\bbb_i}{K}$. Either we have $\ti{A}=\ti{K}$ and \mbox{$\chs_{\st{K}}(\omega)=0$} for $\omega\in\ti{A}$ or $\ti{A}\subsetneq\st{K}$ and $\chs_{\st{K}}(\omega)=1$ for $\omega\in\ti{A}$. Moreover for $K_0\in \GC{j}{\bbb_i}{\Omega}$ and $K_0\neq K$ we know that either $K\subseteq\st{K}_0$ or $K\cap\st{K}_0=\emptyset$. Therefore function $g_j=\sum_{\ti{B}\in \GC{j}{\bbb_i}{\Omega}}  x_{\ti{B}}\Pp(\ti{B})\chs_{\st{B}}$ is constant on $\ti{A}$.}
    \label{fig:constsupp}
\end{figure}
Next put 
\[
\DD{j}=\Est{\GC{j-1}{\bbb_i}{\Omega}}
\]
for $j\geq 2$ and $\DD{1}=\Omega$.
We claim that
\begin{equation}\label{eq: geodecdiam}
   \int_{\DD{j+k+1}}\nx{g_j}^p\dP\leq 2^{-k+1} \eLpx{g_j}^p.   
\end{equation}
In order to see that \eqref{eq: geodecdiam} holds true it suffices to observe that for $\ti{B}\in\GC{j+1}{\bbb_i}{\Omega}$ we have
\[
\sum_{\ti{A}\in\GC{k}{\bbb_i}{\ti{B}}} \Pp(\ti{A}) \stackrel{\eqref{geomet}}{\leq} 2^{-k} \Pp(\ti{B}).
\]
(For a more detailed verification of \eqref{eq: geodecdiam} see \eqref{smallgeoint} above.)
Therefore applying Lemma \ref{estDD} with
\[
a_k=\left\{
\begin{array}{cc} 
     2^{\frac{2-k}{p}} & k\geq 3 \\
     1& k\in\{1,2\}
\end{array}
\right.
\]
gives
\[
\begin{split}
\eLpx{f_i}^p&\leq \left(1+\frac{1}{1-2^{-\frac{1}{p}}}\right)\sum_{j\in\NN}\eLpx{g_j}^{p}
\\&\stackrel{\eqref{highordcol}}{\leq}\left(1+\frac{1}{1-2^{-\frac{1}{p}}}\right)\stlemp\sum_{j\in\NN}\sum_{\ti{A}\in \GC{j}{\bbb_i}{\Omega}}  \nx{x_{\ti{A}}}^p\eLp{\hai{A}}^p
\\&=\left(1+\frac{1}{1-2^{-\frac{1}{p}}}\right)\stlemp\sum_{\ti{A}\in \bbb_i}  \nx{x_{\ti{A}}}^p\eLp{\hai{A}}^p.
\end{split}
\]

\end{proof}
\subsection{Proof of Theorem \ref{thmfincarl}}
Now we have all of the tools for the proof of Theorem \ref{thmfincarl}. 
\begin{proof}[Proof of Theorem \ref{thmfincarl}] Recall that the collection $\Ec$ was defined in \eqref{defcole} and $\{\ha{A}\}_{A\in\Ec}$ in \eqref{defbasis}. Let $f\in L^p(\Omega,X)$. For simplicity of notation we assume that $\EE f=0$. Note that there exists a unique sequence $\{x_A\}_{A\in\Ec}$ such that
\[
f(\omega)=\sum_{A\in\Ec} x_A \hax{A}{\omega}
\]
 where the series on the right hand side converges in $L^p(\Omega,X)$.
Clearly for $n\in\NN$ we have
\[
\mardifx{f}{n}{\omega}=\sum_{A\in \Ec \cap \AAA_n} x_{A} \hax{A}{\omega}.
\]
Next observe that there exists a unique sequence $\{y_A\}_{A\in\AAA}$ in $X$ such that for any $n\in\NN$ the following identity holds
\begin{equation}\label{defygrekow}
\sum_{A\in \Ec \cap \AAA_n} x_{A} \hax{A}{\omega}=\sum_{A\in \AAA_{n}} y_{A} \chs_{A}(\omega). \end{equation}
Recall that in \eqref{collectb} we defined $\bbb=\{\ti{A}: A\in\AAA \}$, where $\ti{A}$ is specified in \eqref{deftiAsun}. Recall that given $A\in\AAA$ we defined $\st{A}=A\setminus\ti{A}.$
For any set $\ti{A}\in\bbb$ we define
\begin{equation}\label{defytilde}
z_{\ti{A}}=\left\{\begin{array}{cc}
  \frac{y_{\st{A}}}{\Pp(\ti{A})}   & \mbox{when } \Pp(\ti{A})\neq 0 \\
    0 & \mbox{when }\Pp(\ti{A})= 0
\end{array}\right.     
\end{equation}
  Recall that we put $\ccc=\{\st{A}: A\in\AAA\}$. We observe that $\ccc\subset \AAA$, $\Ec\setminus\ccc=\emptyset$ and $\AAA=\Ec\cup\,\ccc$. Accordingly we decompose $f$ as $f= g+b$, where
\[
g(\omega)= \sum_{A\in \Ec} y_{A} \chs_{A}(\omega)\qquad\mbox{and}\qquad b(\omega)=\sum_{\st{A}\in\ccc} y_{\st{A}} \chs_{\st{A}}(\omega).
\]
Clearly when $\Pp(\ti{A})= 0$ we have $y_{\st{A}}=0=\Pp(\ti{A})z_{\ti{A}}$. Thus
\[
b(\omega)=\sum_{\ti{A}\in\bbb} z_{\ti{A}} \chs_{\st{A}}(\omega)\Pp(\ti{A}).
\]
We have an estimate for the norm of $f$:
\begin{equation}\label{eq: fprzezgb}
 \eLpx{f}^p\leq 2^{p-1}\left(  \eLpx{g}^p+\eLpx{b}^p\right).   
\end{equation}
By Lemma \ref{conscol}
\begin{equation}\label{estpartcol}
 \begin{split}
\eLpx{g}^{p} &\leq \stlemcrlp{\Ec} \sum_{A\in \Ec} \nx{y_{A}}^p\Pp(A).
\end{split}   
\end{equation}
We use Lemma \ref{estcolb} to obtain that for the function $b$:
\begin{equation*}
\eLpx{b}^p\leq \tilde{C}_p^p\sum_{\ti{A}\in \bbb} \nx{z_{\ti{A}}}^p\eLp{\hai{A}}^p .    
\end{equation*}
where $\hai{A}$ is defined in \eqref{defrozhaa}. Recall that $\eLp{\hai{A}}^p= \Pp(\ti{A})^p\Pp(\st{A})+\eLpx{\Pp(\st{A})\chs_{\ti{A}}}^p$. Thus we have the following estimate
\begin{equation}\label{eq: estbfunct}
\eLpx{b}^p\leq \tilde{C}_p^p\left( \sum_{A\in \AAA} \nx{y_{\st{A}}}^p\Pp(\st{A})+\sum_{\ti{A}\in \bbb} \eLpx{z_{\ti{A}}\Pp(\st{A})\chs_{\ti{A}}}^p\right) .    
\end{equation}
Recall that in \eqref{defcole} we observe that $\Ec=\bigcup_{A\in\AAA}\bigcup_{j=2}^{N(A)} A_j$. Combining \eqref{eq: fprzezgb} with \eqref{estpartcol}, \eqref{eq: estbfunct} and taking in to account that
\[
\tilde{C}^p_p>\stlemcrlp{\Ec}
\]
we have
\[
\begin{split}
\frac{\eLpx{f}^p}{2^{p-1} \tilde{C}^p_p} &\leq  \sum_{A\in \AAA}\sum_{j=2}^{\Nay{A}}\nx{y_{A_j}}^p\Pp(A_j) +\sum_{A\in \AAA} \left(\nx{y_{\st{A}}}^p\Pp(\st{A}) +  \eLpx{z_{\ti{A}}\Pp(\st{A})\chs_{\ti{A}}}^p \right).
\end{split}
\]
Since for any $A\in\AAA$ we have $\st{A}=A_1$ we have
\[
\begin{split}
\sum_{j=2}^{\Nay{A}}\nx{y_{A_j}}^p\Pp(A_j) &+ \nx{y_{\st{A}}}^p\Pp(\st{A})=\sum_{j=1}^{\Nay{A}}\nx{y_{A_j}}^p\Pp(A_j).
\end{split}
\]
Therefore
\[
\begin{split}
\sum_{A\in \AAA}\sum_{j=2}^{\Nay{A}}\nx{y_{A_j}}^p\Pp(A_j) &+\sum_{A\in \AAA} \nx{y_{\st{A}}}^p\Pp(\st{A})=\sum_{A\in\AAA}\nx{y_{A}}^p\Pp(A)
\\&=\sum_{n=1}^{\infty}\sum_{A\in\AAA_n}\nx{y_{A}}^p\Pp(A).
\end{split}
\]
Since we are interested in martingale difference we point out that (see \eqref{defygrekow} and equation before it): 
\begin{equation}\label{eq: ymardif}
 \sum_{n=1}^{\infty}\sum_{A\in\AAA_n}\nx{y_{A}}^p\Pp(A)= \sum_{n=1}^{\infty}\eLpx{\mardif{f}{n}}^p .   
\end{equation}
Now we estimate the remaining part in the series. Recall that $\sum_{j=1}^{\Nay{A}} y_{A_j}\Pp(A_j)=0$. For $z_{\ti{A}}\neq 0$ we have
\[
z_{\ti{A}}\Pp(\st{A})=-\frac{\sum_{j=2}^{\Nay{A}} y_{A_j}\Pp(A_j)}{\Pp(\ti{A})}.
\]
Thus
\[
\begin{split}
 \eLpx{z_{\ti{A}}\Pp(\st{A})\chs_{\ti{A}}}^p&=\Pp(\ti{A})\nx{\frac{\sum_{j=2}^{\Nay{A}} y_{A_j}\Pp(A_j)}{\Pp(\ti{A})}}^p.
\end{split}
\]
By triangle inequality we have
\[
\begin{split}
 \eLpx{z_{\ti{A}}\Pp(\st{A})\chs_{\ti{A}}}^p &\leq \Pp(\ti{A})\left(\sum_{j=2}^{\Nay{A}} \nx{y_{A_j}}\frac{\Pp(A_j)}{\Pp(\ti{A})}\right)^p  
\\&= \left(\sum_{j=2}^{\Nay{A}} \nx{y_{A_j}}\Pp(A_j)^{\frac{1}{p}}\left(\frac{\Pp(A_j)}{\Pp(\ti{A})}\right)^{1-\frac{1}{p}}\right)^p  .
\end{split}
\]
Applying H\"older's inequality gives
\[
\begin{split}
 \eLpx{z_{\ti{A}}\Pp(\st{A})\chs_{\ti{A}}}^p&\leq\left(\sum_{j=2}^{\Nay{A}} \nx{y_{A_j}}^p\Pp(A_j)\right)\left(\sum_{j=2}^{\Nay{A}}\frac{\Pp(A_j)}{\Pp(\ti{A})}\right)^{p-1} 
\\ &\stackrel{\eqref{deftiAsun}}{=}\sum_{j=2}^{\Nay{A}} \nx{y_{A_j}}^p\Pp(A_j).
 \end{split}
\]
Summing the inequalities over $\ti{A}\in \bbb$ we get
\[
\sum_{\ti{A}\in \bbb} \eLpx{y_{\ti{A}}\Pp(\st{A})\chs_{\ti{A}}}^p\leq  \sum_{n=1}^{\infty}\sum_{A\in\AAA_n}\nx{y_{A}}^p\Pp(A).
\]
Therefore
\[
\eLpx{f}^p\leq 2^{p} \tilde{C}^p_p \sum_{n=1}^{\infty}\sum_{A\in\AAA_n}\nx{y_{A}}^p\Pp(A).
\]
Last equation follows from the pairwise disjointedness of sets from $\AAA_n$. By \eqref{eq: ymardif}
\[
\eLpx{f}^p\leq 2 \tilde{C}_p \left(\sum_{n=1}^{\infty}\eLpx{\mardif{f}{n}}^p   \right)^{\frac{1}{p}}
.\]
\end{proof}

\section{Closing remarks on UMD spaces}
In this closing section we consider only Banach spaces with the UMD property. We show that for the class of UMD-Banach spaces Maurey's isomorphism may be employed to obtain a proof of Theorem \ref{main} by direct reduction to the results in \cite{Geiss2008}. This is outlined here in four separate steps.

\paragraph{Review of \cite{Geiss2008}:} Let $\Rc$ be a nested collection of subsets of a probability space $(\Omega,\; \FF,\;\Pp)$. We assume that $\Rc$ satisfies the following additional property. 
\[
\mbox{If}\quad I,J\in\Rc \quad\mbox{and}\quad I\subset J,\; I\neq J\quad\mbox{then}\quad\Pp(I)\leq \frac{\Pp(J)}{2}. 
\]
Let $\{r_I:[0,1]\rightarrow \{-1,1\}\}_{I\in\Rc}$ be an enumeration of an independent Rademacher system; and define
\[
d_I(\omega,t)=\chs_{I}(\omega) r_I(t),\quad \omega\in\Omega,\quad t\in [0,1]
\]
or shortly $d_I=\chs_{I}\otimes\; r_I$.

Let $X$ be a Banach space. We define $X_p(\Rc)$ to be the closure of 
\[
\operatorname{span}\{x_I d_I\;:\; x_I\in X,\;I\in \Rc\}
\]
in $L^{p}(\Omega \times [0,1],\,\FF\otimes \dyadsigma{},\,\Pp\otimes\lambda\;, X)$.  Let the norm on $X_p(\Rc)$ be the one induced by $L^{p}(\Omega \times [0,1], X)$.

The proof given in \cite{Geiss2008} yields that the following dichotomy holds for the class of $X_p(\Rc)$ spaces:
\begin{itemize}
    \item If $\crl{\Rc}<\infty$, then there exists $C(\crl{\Rc},p)>0$ such that
\begin{equation}\label{eq: GMinfcarl}
\|\sum_{I\in \Rc} x_I d_I\|^p_{X_p(\Rc)}\leq C(\crl{\Rc},p) \sum_{I\in\Rc} \nx{x_I}^p |I|,    
\end{equation}
for any Banach space $X$, where $x_I\in X$, $I\in\Rc$ and the series $\sum_{I\in \Rc} x_I d_I$ converges in $L^{p}(\Omega\times [0,1], X)$. 
\item If $\crl{\Rc}=\infty$, then the validity of the estimate
\[
\|\sum_{I\in \Rc} x_I d_I\|^p_{X_p(\Rc)}\leq K^p \sum_{I\in\Rc} \nx{x_I}^p |I|,\qquad x_I\in X
\]
implies that $X$ is of Haar type $p$. 
\end{itemize}

\paragraph{Review of Maurey's isomorphism:} Let $(\Omega,\; \FF,\;(\FF_n),\;\Pp)$ be a filtered probability space. Let $\Ec\subset\FF$ be defined by \eqref{defcole}. As shown in \cite{Mueller1991} (see also \cite[Section 4.1]{Mueller2005}), Maurey's proof of \cite[Prop 4.6, Lemma 4.10]{Maurey2} can be extended to the vector valued case as follows: if X satisfies the UMD property then there exists an isomorphism $T: L^{p}(\Omega,\,\FF,\,\Pp\,; X)\rightarrow X_p(\Ec)$ such that for $f\in L^{p}(\Omega,\,\FF,\,\Pp\,; X)$,
\[
T\mardifold{f}{n}\in \overline{\operatorname{span}}\{x_I d_I\;:\; I\in \GC{n}{\Ec}{\Omega}, \; x_I\in X \},
\]
where $\mardifold{f}{n}=\mardif{f}{n}$. Specifically if $f\in L^{p}(\Omega,\,\FF,\,\Pp\,; X)$  and
\[
Tf= \sum_{I\in\Ec} x_I d_I
\]
then 
\[
T\mardifold{f}{n}= \sum_{I\in \GC{n}{\Ec}{\Omega}}x_I d_I.
\]Conversely if
\[
f= T^{-1}\left( \sum_{I\in\Ec} x_I d_I\right)
\]
then 
\[
\mardifold{f}{n}= T^{-1}\bigg(\sum_{I\in \GC{n}{\Ec}{\Omega}}x_I d_I\bigg).
\]
Moreover
\[
\|T\|\|T^{-1}\|\leq C(p, UMD(X,p)).
\]

\paragraph{Finite Carleson constant:} Fix a filtered probability space $(\Omega,\,\FF,\,(\FF_n),\,\Pp)$ such that $\crl{\Ec}<\infty$. Fix $1<p<\infty$. Choose $f\in L^p (\Omega,\,\FF,\,\Pp\,; X)$ and let $x_I$ be such that
\[
Tf= \sum_{I\in\Ec} x_I d_I.
\]
By \eqref{eq: GMinfcarl} we have 
\[
\|\sum_{I\in \Ec} x_I d_I\|^p_{X_p(\Ec)}\leq C(\crl{\Ec},p) \sum_{I\in\Ec} \nx{x_I}^p |I|.
\] 
Hence,
\[
\begin{split}
\eLpx{f}^p&\leq \|T^{-1}\|^p\|Tf\|_{X_p(\Ec)}^p=\|T^{-1}\|^p\bigg\|\sum_{I\in\Ec} x_I d_I\bigg\|_{X_p(\Ec)}^p
\\&\leq C(\crl{\Ec},p)\|T^{-1}\|^p \sum_{I\in \Ec}\nx{x_I}^p |I|
\\&= C(\crl{\Ec},p)\|T^{-1}\|^p\sum_{n} \sum_{I\in \GC{n}{\Ec}{\Omega}}\nx{x_I}^p |I|
\\&= C(\crl{\Ec},p)\|T^{-1}\|^p\sum_{n} \bigg\|\sum_{I\in \GC{n}{\Ec}{\Omega}}x_I d_I\bigg\|_{X_p(\Ec)}^p
\\&= C(\crl{\Ec},p)\|T^{-1}\|^p\sum_{n} \|T\mardifold{f}{n}\|_{X_p(\Ec)}^p
\\&\leq C(\crl{\Ec},p)\|T^{-1}\|^p \|T\|^p\sum_{n} \eLpx{\mardifold{f}{n}}^p.
\end{split}
\]

\paragraph{Infinite Carleson constant:} Conversely assume that $\crl{\Ec}=\infty$. Let $1<p\leq 2$. Assume that for every $f\in  L^p (\Omega,\,\FF,\,\Pp\,; X)$ we have
\[
\eLpx{f}^p\leq K^p \sum_{n} \eLpx{\mardifold{f}{n}}^p.
\]
Fix $g\in X_p(\Ec)$ with expansion
\[
g=\sum_{I\in\Ec} x_I d_I.
\]
Apply the Maurey isomorphism to $g$ and put $f=T^{-1} g.$
\[
\begin{split}
\bigg\|\sum_{I\in\Ec} x_I d_I\bigg\|_{X_p(\Ec)}^p&=\|Tf\|^p_{X_p(\Ec)} \leq \|T\|^p\eLpx{f}^p
\\&\leq K^p\|T\|^p\sum_{n} \eLpx{\mardifold{f}{n}}^p
\\&= K^p\|T\|^p\sum_{n} \eLpxb{T^{-1}\sum_{I\in \GC{n}{\Ec}{\Omega}}x_I d_I}^p
\\&\leq K^p\|T\|^p\|T^{-1}\|^p\sum_{n} \bigg\|\sum_{I\in \GC{n}{\Ec}{\Omega}}x_I d_I\bigg\|_{X_p(\Ec)}^p
\\&= K^p\|T\|^p\|T^{-1}\|^p \sum_{I\in \Ec}\nx{x_I}^p |I|.
\end{split}
\]
In view of \cite{Geiss2008} $X$ satisfies $HT_p$ and $MT_p$.\qed \vspace{0.5cm}

{\bf Acknowledgements:} We are indebted to Maciej Rzeszut for very instructive and informative conversations during preparation of this paper and for his permission to include the example presented in the Remark \ref{MREX}.

 We are most grateful to the anonymous reviewer. Their extremely thoughtful feedback have significantly improved our paper. We truly appreciate their highly valuable contribution.
\bibliographystyle{abbrv}
\bibliography{martingale}
Institute of Analysis,\\ Johannes Kepler University Linz,\\ Austria, 4040 Linz,\\ Altenberger Strasse 69\\
E-mail address: paul.mueller@jku.at\\\textcolor{white}{,}\\
Institute of Analysis, \\Johannes Kepler University Linz,\\ Austria, 4040 Linz,\\ Altenberger Strasse 69\\
and \\
Institute of Mathematics,\\ University of Warsaw,\\ Poland, 02-097 Warszawa,\\ ul. Banacha 2\\
E-mail address: krystian.kazaniecki@jku.at
\end{document}